\documentclass[12pt]{amsart}
\usepackage[utf8]{inputenc}

\usepackage{subfiles}
\usepackage{hyperref}
\usepackage{cleveref}
\usepackage{amsmath}
\usepackage{amssymb}
\usepackage{amsthm}
\usepackage{mathtools}
\usepackage{thmtools}
\usepackage{thm-restate}
\usepackage[margin=1in]{geometry}
\usepackage{cleveref}
\usepackage{ytableau}
\usepackage{comment}
\usepackage{soul}

\usepackage[maxnames=5]{biblatex}
\renewbibmacro{in:}{}
\DeclareFieldFormat[article]{citetitle}{#1}
\DeclareFieldFormat[article]{title}{#1}

\numberwithin{equation}{section}

\theoremstyle{definition}
\newtheorem{theorem}{Theorem}[section]
\newtheorem*{theorem*}{Theorem}
\newtheorem{example}[theorem]{Example}
\newtheorem*{example*}{Example}
\newtheorem{lemma}[theorem]{Lemma}
\newtheorem*{lemma*}{Lemma}
\newtheorem{corollary}[theorem]{Corollary}
\newtheorem*{corollary*}{Corollary}
\newtheorem{definition}[theorem]{Definition}
\newtheorem*{definition*}{Definition}

\newtheorem*{proposition*}{Proposition}

\newtheorem*{remark*}{Remark}
\newtheorem{conjecture}[theorem]{Conjecture}

\usepackage{ mathrsfs }

\title{The Image of the Pop Operator on Various Lattices}
\author{Yunseo Choi}\address{\textsc{Y. Choi}, Harvard University,
    Cambridge, MA, 02138} \email{ychoi@college.harvard.edu}
\author{Nathan Sun}\address{\textsc{N. Sun}, Harvard University,
    Cambridge, MA, 02138} \email{nsun@college.harvard.edu}

\begin{document}

\maketitle

\begin{abstract}
Extending the classical pop-stack sorting map on the lattice given by the right weak order on $S_n$, Defant defined, for any lattice $M$, a map  $\mathsf{Pop}_{M}: M \to M$ that sends an element $x\in M$ to the meet of $x$ and the elements covered by $x$. In parallel with the line of studies on the image of the classical pop-stack sorting map, we study $\mathsf{Pop}_{M}(M)$ when $M$ is the weak order of type $B_n$, the Tamari lattice of type $B_n$, the lattice of order ideals of the root poset of type $A_n$, and the lattice of order ideals of the root poset of type $B_n$. In particular, we settle four conjectures proposed by Defant and Williams on the generating function 
\begin{equation*}
    \mathsf{Pop}(M; q) = \sum_{b \in \mathsf{Pop}_{M}(M)} q^{|\mathscr{U}_{M}(b)|},
\end{equation*}
where $\mathscr{U}_{M}(b)$ is the set of elements of $M$ that cover $b$. 
\end{abstract}

\section{Introduction} 
\label{intro}
In 1982, Ungar \cite{ungar19822n} defined the pop-stack sorting map to act on the permutations in the symmetric group $S_n$ by reversing their descending runs. Since then, the pop-stack sorting map has been thoroughly studied \cite{pop1, pop2, pop3, claesson2019enumerating, pudwell2018two}. 

More recently in 2022, Defant \cite{def} generalized the classical pop-stack sorting map so that it is defined on any lattice $M$. Defant's map $\mathsf{Pop}_{M}$ sends each element $x\in M$ to the meet of $x$ and the set of elements covered by $x$. In particular, when $M$ is the right weak order on $S_n$, Defant's $\mathsf{Pop}_{M}$ coincides with Ungar's pop-stack sorting map. In \cite{def}, Defant studied $\mathsf{Pop}_M$ when $M$ is the $\nu$-Tamari lattice; in \cite{popstacksorting}, he studied the case when $M$ is the weak order on a finite Coxeter group.

In \cite{conjecture}, Defant and Williams proved a theorem that related the image of $\mathsf{Pop}_M$ with the image of the \emph{dual pop-stack sorting operator}, the \emph{rowmotion} operator, and \emph{independent dominating sets}. Motivated by their theorem, Defant and Williams defined the following generating function that refines the count of elements of $M$ that lie in the image of $\mathsf{Pop}_{M}$: 
\begin{equation*}
    \mathsf{Pop}(M; q) = \sum_{b \in \mathsf{Pop}_{M}(M)} q^{|\mathscr{U}_{M}(b)|},
\end{equation*}
where $|\mathscr{U}_{M}(b)|$ is the number of elements of $M$ that cover $b$. 

The analog of $\mathsf{Pop}(M; q)$ for the pop-stack sorting map was previously studied \cite{pop1, cgp}. Notably, in Asinowski, Banderier, Billey, Hackl, and Linusson \cite{pop1}, it was proved that $$[q^{n-2}] \mathsf{Pop}(\mathrm{Weak}(A_{n-1}); q) = 2^{n} - 2n,$$ where $\mathrm{Weak}(A_{n-1})$ is the weak order on the symmetric group $S_n$.  In addition, when $M = J(\Phi^{+}_{A_n})$ is the lattice of order ideals of the root poset of type $A_n$, Sapounakis, Tasoulas, and Tsikouras \cite{greek} proved that 
\begin{equation*}
    \mathsf{Pop}(J(\Phi^{+}_{A_n}); 1) = \sum_{k =0}^{n} \frac{1}{k+1} \sum_{j=0}^{n-k+1} \binom{k+1}{j-1} \binom{k+1}{j} \binom{n-j+1}{n-k-j+1}.  
\end{equation*}
More recently, Hong \cite{carina} settled a conjecture of Defant and Williams \cite{conjecture} by proving that 
\begin{equation*}
    \mathsf{Pop}(\mathrm{Tam}(A_n);q) = \sum_{k=0}^{n} \frac{1}{k+1} \binom{2k}{k} \binom{n}{2k} q^{n-k},
\end{equation*}
where $\mathrm{Tam}(A_n)$ is the Tamari lattice of type $A_n$. Specializing $q=1$ in this result tells us that the size of the image of $\mathsf{Pop}_{\mathrm{Tam}(A_n)}$ is the $n$th Motzkin number. 

In this paper, we settle four conjectures about the generating functions $\mathsf{Pop}(M; q)$ proposed by Defant and Williams \cite{conjecture}. In what follows, $\mathrm{Weak}(B_n)$ is the weak order of type $B_n$, $\mathrm{Tam}(B_n)$ is the Tamari lattice of type $B_n$, and $J(\Phi_{A_n}^+)$ and $J(\Phi_{B_n}^+)$ are the lattices of order ideals of the root posets of types $A_n$ and $B_n$, respectively.


\begin{theorem}
\label{weak}
For all $n$, we have that 
\begin{equation*}
    [q^{n-1}] \mathsf{Pop}\mathrm{(Weak}(B_n); q) = 3^{n}-2n-1. 
\end{equation*}
\end{theorem}


\begin{theorem}\label{tamaribn} For all $n$, we have that
$$\mathsf{Pop}(\mathrm{Tam}(B_n); q) = \sum_{k =0}^{\lfloor \frac{n+1}{2} \rfloor} \binom{n-1}{k} \binom{n+1-k}{k} q^{n-k} .$$
\end{theorem}


\begin{theorem} \label{jayan}
For all $n$, we have that 
\begin{equation*}
    \mathsf{Pop}(J(\Phi^{+}_{A_n}); q) = \sum_{k =0}^{n} \frac{1}{k+1} \sum_{j=0}^{n-k+1} \binom{k+1}{j-1} \binom{k+1}{j} \binom{n-j+1}{n-k-j+1} q^{k+1}.  
\end{equation*}
\end{theorem}


\begin{theorem} \label{jaybn}
For all $n$, we have that 
\begin{equation*}
    \mathsf{Pop}(J(\Phi^{+}_{B_n}); q) = \sum_{k = 0}^{n} \sum_{j=1}^{k} \binom{2j}{j} \binom{k+j}{k-j} \binom{n-k+j-1}{n-k} (-1)^{k-j} q^{k}.  
\end{equation*}
\end{theorem}

The rest of the paper is organized as follows. In \Cref{prelim}, we formally define $\mathsf{Pop}_{M}$ on lattices $M$ and describe the action of $\mathsf{Pop}_{M}$ when $M$ is $\mathrm{Weak}(B_n)$, $\mathrm{Tam}(B_n)$, $J(\Phi^{+}_{A_n})$, and $J(\Phi^{+}_{B_n})$. In \Cref{proofs}, we present the proofs of our main results. In \Cref{futuredirections}, we suggest potential future directions of study. 

\section{Preliminaries}
\label{prelim}
\subsection{Definitions}

Suppose that $x = x_{1} x_{2} \cdots x_{n},$ where the natural numbers $x_{i}$ are distinct. The \emph{length} of $x$ is $\mathrm{len}(x) = n$. In addition, the \emph{reduction} of $x$ is $\mathrm{red}(x) = x_1'\cdots x_n' \in S_{n}$, where $x_i' = i$ if $x_{i}$ is the $i^{\text{th}}$ smallest number from the set $\{x_{1}, x_{2}, \ldots, x_{n}\}$. For example, if $x = 8512$, then $\mathrm{red}(x) = 4312$. Furthermore, if $x = x_{1} x_{2} \cdots x_{n} \in S_{n}$, then let $\mathrm{ind}_{x}(i)$ be such that $x_{\mathrm{ind_{x}(i)}} = i$. For example, if $x = 4312$, then $\mathrm{ind}_{x}(3) = 2$. 

A \textit{descent} of a permutation $x = x_{1} x_{2} \cdots x_{n} \in S_{n}$ is an adjacent pair of numbers $x_{i}$ and $x_{i+1}$ such that $x_{i} > x_{i+1}$. For $x=51763284$, the descents of $x$ are $(5, 1)$, $(7,6)$, $(6,3)$, $(3,2)$, and $(8,4)$. We let $|\mathrm{des}(x)|$ represent the number of descents in $x$. Continuing with our example of $x = 51763284$, we have that $|\mathrm{des}(x)| = 5$. A \textit{descending run} is a maximal decreasing substring of $x$. For $x=51763284$, the descending runs are $51$, $7632$, and $84$. Next, let $\mathrm{rev}(x)$ denote the permutation that results after reversing each descending run of $x$. So, if $x=51763284$, then $\mathrm{rev}(x) = 15236748$. Furthermore, for $x = x_{1} x_{2} \cdots x_{n}$, let $x+k :=(x_{1} +k)(x_{2}+k) \cdots (x_{n}+k)$. For example, if $x = 2143$, then $x+2 =4365$. Additionally, let $\cdot$ denote concatenation.

We say that the permutation $x$ \emph{contains} the pattern $y$ if there exists a sequence of indices $c(1) < \dots < c(k)$ such that $x_{c(1)} \cdots x_{c(k)}$ is order-isomorphic to $y$. Equivalently, we say that $x_{c(1)} \cdots x_{c(k)}$ is a \emph{$y$ pattern of $x$}. For example, $x_1x_3x_5 = 513$ is a $312$ pattern of $x = 52143$. Furthermore, we say that the permutation $x$ contains the \emph{vincular} pattern $y$ if there exists indices $c(1), \dots , c(k)$ such that $x_{c(1)} \cdots x_{c(k)}$ is a $y$ pattern of $x$ and certain indices $c(1), c(2), \ldots, c(n)$ are consecutive. The entries of $y$ that correspond to the indices that must be consecutive are overlined in $y$. Equivalently, we say that $x_{c(1)} \cdots x_{c(k)}$ is a vincular $y$ pattern of $x$. For example, $x_1x_2x_5 = 523$ in $x = 52143$ is a vincular $\overline{31}2$ pattern of $x$, while $x_1x_3x_5=513$ is not, because $x_1$ and $x_3$ are not consecutive in $x$. We say a permutation \emph{avoids} a (vincular) pattern if it does not contain it. For instance, the permutation $x = 23145$ avoids the pattern $312$ pattern as well as the vincular pattern $\overline{31}2$. 

Lastly, we say the subsequence $x_{i(1)} x_{i(2)} x_{i(3)}$ of $x = x_{1} x_{2} \cdots x_{2n}$ is a \textit{$312^{*}$ pattern} of $x$ if  $\mathrm{red}(x_{i(1)} x_{i(2)} x_{i(3)}) = 312$ and $x_{i(3)} \geq n+1$. For example, there is no $312^{*}$ pattern in $x = 3142$, even though $312$ is a $312$ pattern of $x$. We say that the subsequence $x_{i(1)}x_{i(2)}x_{i(3)}$ of $x$ is a \textit{big $312^{*}$ pattern} of $x$ if $x_{i(1)} x_{i(2)} x_{i(3)}$ is a $312^{*}$ pattern and $x_{i(2)} \geq n+1$. Conversely, we say that the subsequence $x_{i(1)}x_{i(2)}x_{i(3)}$ of $x$ is a \textit{small $312^{*}$ pattern} of $x$ if $x_{i(1)} x_{i(2)} x_{i(3)}$ is a $312^{*}$ pattern and $x_{i(2)} \leq n$. In a related direction, we define $x_{i(1)}x_{i(2)}x_{i(3)}$ to be a $213^{*}$ pattern if $\mathrm{red}(x_{i(1)} x_{i(2)} x_{i(3)}) = 213$ and $x_{i(1)} \leq n$.

\subsection{Lattice basics and the \texorpdfstring{$\mathsf{Pop}$}{Pop} operator} 

A poset $M$ is a \textit{lattice} if any two elements $x, y \in M$ have a greatest lower bound, which we call their \textit{meet} and denote as $x \wedge y$, and a least upper bound, which we call their \textit{join} and denote as $x\vee y$. We let $\preceq$ denote the partial order on $M$. We say that $L \subset M$ is a \emph{sublattice} of $M$ if $L$ and $M$ have the same meet $\wedge$ and join $\vee$ operations and for any $x, y \in L$, we have $x \wedge y, x \vee y \in L$. 

For $x, y \in M$, we say that $y$ covers $x$---which we denote by $x \lessdot y$---if $x \prec y$, and there is no $z \in M$ such that $x \prec z \prec y$. For $x \in M$, let $\mathscr{U}_{M}(x)$ denote the set of elements of $M$ that cover $x$. In a related direction, for $x \in M$, let $\mathscr{D}_{M}(x)$ denote the set of elements of $M$ that $x$ covers. 

An equivalence relation $\equiv_{L}$ on $M$ is a \textit{lattice congruence} if for all $x_{1}, x_{2}, y_{1}, y_{2} \in M$ such that $x_{1} \equiv_{L} x_{2}$ and $y_{1} \equiv_{L} y_{2}$, we have that $x_{1} \wedge y_{1} \equiv_{L} x_{2} \wedge y_{2}$ and $x_{1} \vee x_{2} \equiv_{L} y_{1} \vee y_{2}$. It is well known that if $M$ is a finite lattice and $\equiv_{L}$ is a lattice congruence of $M$, then every congruence class of $\equiv_L$ is an interval. Therefore, we can define $\pi_{L_\downarrow}(x)$ to be the unique minimal element in the congruence class of $x$ defined by $\equiv_{L}$. 

We now define the action of the \textit{pop-stack sorting operator} $\mathsf{Pop}^{\downarrow}_{M}: M \to M$. 

\begin{definition}[Defant, \cite{def}]
Given a lattice $M$ and an element $x \in M$, we let 
$$\mathsf{Pop}^{\downarrow}_M(x):=\bigwedge(\{y\in M: y\lessdot x\}\cup \{x\}).$$
\end{definition}

For the rest of the paper, we use $\mathsf{Pop}_{M}$ to mean $\mathsf{Pop}^{\downarrow}_{M}$. In a related direction, we define the action of the \textit{dual pop-stack sorting operator} $\mathsf{Pop}^{\uparrow}_{M}: M \to M$. 

\begin{definition}[Defant and Williams, \cite{conjecture}]
    Given a lattice $M$ and an element $x \in M$, let 
    \begin{equation*}
        \mathsf{Pop}^{\uparrow}_{M}(x):=\bigvee(\{y \in M: x \lessdot y\} \cup \{x\}).  
    \end{equation*}
\end{definition}

Next we define the generating function that refines the count of the elements of $M$ that belong to $\mathsf{Pop}_{M}(M)$. 

\begin{definition}[Defant and Williams, \cite{conjecture}] For a lattice $M$, let 
        \begin{equation*} 
            \mathsf{Pop}(M; q) := \sum_{x \in \mathsf{Pop}_{M}(M)} q^{|\mathscr{U}_{M}(x)|}.
        \end{equation*}
\end{definition}

We end this subsection by citing the following theorem by Defant and Williams \cite{conjecture} that shows how to write $\mathsf{Pop}(M; q)$ in terms of the dual pop-stack sorting operator. 

\begin{theorem}[Defant and Williams \cite{conjecture}] \label{dual}
    For every lattice $M$, we have that 
    \begin{equation*}
        \mathsf{Pop}(M; q) = \sum_{x \in \mathsf{Pop}_M(M)} q^{|\mathscr{U}_{M}(x)|} = \sum_{x \in \mathsf{Pop}^{\uparrow}_{M}(M)} q^{|\mathscr{D}_{M}(x)|}. 
    \end{equation*}
\end{theorem}

\subsection{\texorpdfstring{$\mathrm{Weak}(B_n)$}{Weak({Bn})}, \texorpdfstring{$\mathrm{Tam}(A_n)$}{Tam(An)}, 
\texorpdfstring{$\mathrm{Tam}(B_n)$}{Tam(Bn)}, \texorpdfstring{$J(\Phi^{+}_{A_n})$}{J(Phi+{An})}, and \texorpdfstring{$J(\Phi^{+}_{B_n})$}{J(Phi+{Bn})}}

Here, we define the relevant lattices and describe the action of $\mathsf{Pop}_{M}$ on them. First, let $\mathrm{Weak}(A_{n-1})$ be the right weak order on $S_n$. In this lattice, $x$ covers $y$ if $y$ is obtained by reversing a descent of $x$. In this setting, Defant's pop-stack sorting map \cite{def} agrees with Ungar's pop-stack sorting map; it acts on a permutation in $S_n$ by reversing each descending run. 

Next, we define $\mathrm{Weak}(B_n)$. Let $B_n$ be the $n$-th hyperoctahedral group, the subgroup of $S_{2n}$ consisting of permutations $x=x_1 x_2 \cdots x_{2n}$ that satisfy $x_{i} + x_{2n+1-i} = 2n+1$ for all $1 \leq i \leq 2n$. Then $\mathrm{Weak}(B_{n})$ is the sublattice of the right weak order on $S_{2n}$ given by the elements of $B_n$. It is well known that $|\mathscr{U}_{\mathrm{Weak}(B_n)}(x)| =  \left| \{i \in [1, 2, \ldots, n] : x_{i} < x_{i+1} \} \right|$ for $x \in \mathrm{Weak}(B_{n})$. The operator $\mathsf{Pop}_{\mathrm{Weak}(B_n)}$ acts on the elements of $\mathrm{Weak}(B_n)$ by reversing its descending runs \cite{defstacksorting} (i.e., $\mathsf{Pop}_{\mathrm{Weak}(B_n)}$ is the restriction of $\mathsf{Pop}_{\mathrm{Weak}(A_{2n-1})}$ to the sublattice $\mathrm{Weak}(B_n)$).

We next define $\mathrm{Tam}(A_n)$ and the congruence relation $\equiv_{\mathrm{Tam}(A_n)}$ on $S_{n+1}$. The lattice $\mathrm{Tam}(A_n)$ can be realized as a sublattice of $\mathrm{Weak}(A_n)$, in which each element is $312-$avoiding. It is known that $|\mathscr{U}_{\mathrm{Tam}(A_n)}(\pi)| = |\{i \in [1, 2, \ldots, n] : \pi_{i} < \pi_{i+1} \}|$ \cite{def}.  We say that $x, y \in S_{n+1}$ are $\mathrm{Tam}(A_n)$-\textit{adjacent} and write that $x \lhd_{\mathrm{Tam}(A_n)} y$ if there exists a sequence of numbers $X$, $Y$, and $Z$ and numbers $a < b< c$ such that $x=XcaYbZ$ and $y= XacYbZ$. Then for $x, y \in S_{n+1}$, we say that $x \equiv_{\mathrm{Tam}(A_{n})} y$ if and only if there exists a chain of permutations $x =v_{1}, v_{2}, \ldots, v_{\ell} =y$ such that $v_{i} \lhd_{\mathrm{Tam}(A_{n})} v_{i+1}$ or $v_{i+1} \lhd_{\mathrm{Tam}(A_{n})} v_{i}$ for all $1\leq i \leq \ell-1$. Then the set of elements of $\mathrm{Tam}(A_n)$ coincides with the set of $x \in S_{n+1}$ that can be written as $x = \pi_{\mathrm{Tam}(A_n)_\downarrow}(y)$ for some $y \in S_{n+1}$ \cite{def,cambrianlattices}. 

The action of $\mathrm{Pop}_{\mathrm{Tam}(A_n)}$ can now be described as follows. 

\begin{theorem}{(Defant, \cite{def})} \label{pop_tam_a}
For $x \in \mathrm{Tam}(A_n)$, 
\begin{equation*}
    \mathsf{Pop}_{\mathrm{Tam}(A_n)}(x) = \pi_{\mathrm{Tam}(A_{n})_\downarrow}(\mathrm{rev}(x)). 
\end{equation*}
\end{theorem}

Next, we define $\mathrm{Tam}(B_n)$ and the congruence relation $\equiv_{\mathrm{Tam}(B_n)}$ on $B_{n}$. The lattice $\mathrm{Tam}(B_n)$ can be realized as the sublattice of $\mathrm{Weak}(B_n)$ consisting of the permutations that avoid $312^{*}$ \cite{cambrianlattices}. It is known that $|\mathscr{U}_{\mathrm{Weak}(B_n)}(\pi)| = |\{i \in [1, 2, \ldots, n] : \pi_{i} < \pi_{i+1} \}|$ \cite{cambrianlattices}. We say that $x, y \in B_{n}$ are $\mathrm{Tam}(B_n)$-\textit{adjacent} and write that $x \lhd_{\mathrm{Tam}(B_n)} y$ if there exists a sequence of numbers $X$, $Y$, and $Z$ and numbers $a < b< c$ for which $\mathrm{ind}_{x}(b) \geq \mathrm{ind}_{x}(a)$ if $b \geq n+1$ and $\mathrm{ind}_{x}(a) \geq \mathrm{ind}_{x}(b)$ if $b \leq n$ such that if $a+c \ne 2n+1$, then $x=XcaY(2n+1-c)(2n+1-a)Z$ and $y= XacY(2n+1-a)(2n+1-c)Z$ and if $a+c = 2n+1$, then $x = X ca Y$ and $y =X ac  Y$. Then for $x, y \in B_{n}$, we say that $x \equiv_{\mathrm{Tam}(B_{n})} y$ if there exists a chain of permutations $v_{1}, v_{2}, \ldots, v_{\ell}$ such that $x = v_1$, $y = v_\ell$, and $v_{i} \lhd_{\mathrm{Tam}(B_{n})} v_{i+1}$ or $v_{i+1} \lhd_{\mathrm{Tam}(B_{n})} v_{i}$ for all $1\leq i \leq \ell-1$. Then the set of elements of $\mathrm{Tam}(B_n)$ coincides with the set of $x \in B_{n}$ that can be written as $x = \pi_{\mathrm{Tam}(B_n)_\downarrow}(y)$ for some $y \in B_{n}$ \cite{cambrianlattices}.

Now, the action of $\mathsf{Pop}_{\mathrm{Tam}(B_n)}$ can be described as follows. 

\begin{theorem}[Defant, \cite{def}]  \label{pop_tam_b}
For $x \in \mathrm{Tam}(B_n)$, we have that \begin{equation*}
    \mathsf{Pop}_{\mathrm{Tam}(B_n)}(x) = \pi_{\mathrm{Tam}(B_{n})_\downarrow}(\mathrm{rev}(x)). 
\end{equation*}
\end{theorem}

We finish by defining $J(\Phi^{+}_{A_n})$ and $J(\Phi^{+}_{B_n})$. To do so, we introduce Dyck paths. Dyck paths of semi-length $n$ are paths from $(0, 0)$ to $(2n, 0)$ that do not go below the $x-$axis and consists of steps that are either $(1,1)$ or $(1, -1)$. We call the step $(1,1)$ a $\emph{rise}$ and denote it by $\mathit{r}$, while we call the step $(1, -1)$ a \emph{fall} and denote it by $\mathit{f}$. A \emph{valley} (respectively, \emph{peak}) of a Dyck path is a lattice point along the path that is preceded by a fall (respectively, rise) and succeeded by a rise (respectively, fall). Given any Dyck path $\mu$, we let $s(\mu)$ be its semi-length.

Now when $P$ is a finite poset, we write $J(P)$ for the lattice of order ideals of $P$ ordered by inclusion. We write $\Phi_{A_{n-1}}^+$ and $\Phi_{B_n}^+$ for the root posets of types $A_{n-1}$ and $B_n$, respectively. One can identify the elements of $J(\Phi^{+}_{A_{n-1}})$ with the Dyck paths of semi-length $n$. It follows from the definition of the dual pop operator that  $\mathsf{Pop}^{\uparrow}_{J(\Phi^{+}_{A_{n-1}})}$ acts on a path in $\Phi^{+}_{A_{n-1}}$ by changing each valley into a peak. For $\mu \in \Phi^{+}_{A_{n-1}}$, it is well known that $|\mathscr{D}_{J(\Phi^{+}_{A_{n-1}})}(\mu)|$ is given by the number of peaks in $\mu$. 

The lattice $J(\Phi^{+}_{B_n})$ can now be realized as a sublattice of $J(\Phi^{+}_{A_{2n-1}})$ that consists of Dyck paths of semi-length $2n$ that are symmetric about the line $x = 2n$. It again follows from the definition of a dual pop-stack sorting operator that $\mathsf{Pop}^{\uparrow}_{J(\Phi^{+}_{B_n})}$ acts on a path in $\Phi^{+}_{B_{n}}$ by reversing each valley into a peak. For $\mu \in J(\Phi^{+}_{B_{n}})$, it is well known that $|\mathscr{D}_{J(\Phi^{+}_{B_{n}})}(\mu)|$ is given by the number of peaks that has $x$ coordinate at most $2n$. More generally, for a Dyck path $\mu$ of semi-length $n$ that is symmetric about the line $x = n$, we let $p(\mu)$ count the number of peaks along the path that has $x$ coordinate at most $x = n$.

\section{Proofs of the Main Results}
\label{proofs}
\subsection{\texorpdfstring{$\mathrm{Weak}(B_n)$}{Weak(Bn)}} Here, we settle a conjecture by Defant and Williams \cite{conjecture} on the coefficient of $q^{n-1}$ in $\mathsf{Pop}(\mathrm{Weak}(B_n); q)$. We begin with a note on notations that we carry throughout this subsection. 

Suppose that $x = x_{1} x_{2} \cdots x_{2n} \in \mathrm{Weak}(B_{n})$. Let $\mathrm{asc}_{k}(x)$ be the $k^{\text{th}}$ maximal ascending substring of $x$ from the left. For example, if $x=68245713$, then $\mathrm{asc}_1(x) =68$, $\mathrm{asc}_2(x) =2457$, and $\mathrm{asc}_3(x) = 13$. Similarly, let $\mathrm{asc}_{-k}(x)$ be the $k^{\text{th}}$ maximal ascending substring of $x$ from the right. Continuing with our example of $x = 68245713$, we have that $\mathrm{asc}_{-1}(x) =13$, $\mathrm{asc}_{-2}(x) =2457$, $\mathrm{asc}_{-3}(x) = 68$. Furthermore, let $\ell_{k}(x) = \mathrm{len}(\mathrm{asc}_{k}(x))$ so that in our example of $x =68245713$, we have that $\ell_{1}(x) = \ell_3(x) = 2$ and $\ell_2(x) =  4$. 

In addition, for $x \in \mathrm{Weak}(B_n)$, let $\mathrm{mid}(x)$ be the maximal ascending run in $x$ that contains both $x_{n}$ and $x_{n+1}$. If no ascending run contains both $x_{n}$ and $x_{n+1}$, then $\mathrm{mid}(x)$ is the empty string. For $x = 68245713$, we have that $\mathrm{mid}(x) = 2457$. On the other hand, if $x = 68254713$, then we have that $\mathrm{mid}(x)$ is the empty string. 

We now proceed to establish auxiliary lemmas that will build up to the statement of \Cref{weak}. We first establish a necessary condition for $x \in \mathsf{Pop}_{\mathrm{Weak}(B_n)}(\mathrm{Weak}(B_n))$. 

\begin{lemma} \label{first}
Suppose that $x \in \mathsf{Pop}_{\mathrm{Weak}(B_n)}(\mathrm{Weak}(B_n))$. Then for any $k$, we have
\begin{equation}
    (\mathrm{asc}_k(x))_1 < (\mathrm{asc}_{k+1}(x))_{\ell(k+1)}. 
\end{equation}
\end{lemma}
\begin{proof}
     Suppose that for some $x$ such that $(\mathrm{asc}_{k+1}(x))_{\ell(k+1)} < (\mathrm{asc}_k(x))_1$ for some $k$, there exists $y \in \mathrm{Weak}(B_n)$ such that $\mathsf{Pop}_{\mathrm{Weak}(B_n)}(y) = x$. We first know from the condition $(\mathrm{asc}_{k+1}(x))_{\ell(k+1)} < (\mathrm{asc}_k(x))_1$ that every number that appears in $\mathrm{asc}_k(x)$ must be greater than every number that appears in $\mathrm{asc}_{k+1}(x)$. Now because the numbers that appear in the descending run containing $y_{\mathrm{ind}_{x}((\mathrm{asc}_{k}(x))_{\ell(k)})}$ must be a subset of those that appear in $\mathrm{asc}_{k}(x)$ and the numbers that appear in the descending run containing $y_{\mathrm{ind}_{x}((\mathrm{asc}_{k+1}(x))_{1})}$ must be a subset of those that appear in $\mathrm{asc}_{k+1}(x)$, we have that $y_{\mathrm{ind}_{x}((\mathrm{asc}_{k+1}(x))_{1})} < y_{\mathrm{ind}_{x}((\mathrm{asc}_{k}(x))_{\ell(k)})}$. Thus, $y_{\mathrm{ind}_{x}((\mathrm{asc}_{k}(x))_{\ell(k)})}$ and $y_{\mathrm{ind}_{x}((\mathrm{asc}_{k+1}(x))_{1})}$ must belong to the same descending run of $y$. However, if so, then $(\mathrm{asc}_{k}(x))_{\ell(k)}$ and $(\mathrm{asc}_{k+1}(x))_1$ must belong to the same ascending run of $x$. Therefore, we reach a contradiction as sought. 
\end{proof}

We begin counting permutations $x \in \mathsf{Pop}_{\mathrm{Weak}(B_n)}(\mathrm{Weak}(B_n))$ that satisfy $|\mathscr{U}_{\mathrm{Weak}(B_n)}(x)| = n-1$. We start by identifying a class of permutations that belongs to $\mathsf{Pop}_{\mathrm{Weak}(B_n)}(\mathrm{Weak}(B_n))$. 

\begin{lemma} \label{one_one}
Suppose that $y \in \mathsf{Pop}_{\mathrm{Weak}(B_{n-1})}(\mathrm{Weak}(B_{n-1}))$. Then $x =1 \cdot (y+1) \cdot 2n \in \mathsf{Pop}_{\mathrm{Weak}(B_n)}(\mathrm{Weak}(B_{n}))$. 
\end{lemma}
\begin{proof}
    Let $z \in \mathrm{Weak}(B_{n-1})$ such that $\mathsf{Pop}_{\mathrm{Weak}(B_{n-1})}(z) = y$. Now consider $w=1 \cdot (z+1) \cdot 2n$. Because no other number belongs to the same descending run as $1$ or $2n$ in $w$, we have that $\mathsf{Pop}_{\mathrm{Weak}(B_n)}(w) = \mathsf{Pop}_{\mathrm{Weak}(B_n)}(1 \cdot (z+1) \cdot 2n) = 1 \cdot (\mathsf{Pop}_{\mathrm{Weak}(B_n)}(z)+1) \cdot 2n = 1 \cdot (y+1) \cdot 2n =x$ as sought. 
\end{proof}

Next, we show that another class of permutations belong to $\mathsf{Pop}_{\mathrm{Weak}(B_n)}(\mathrm{Weak}(B_n))$. 

\begin{lemma} \label{one_two}
For each $1 \leq j \leq n-1$, the permutation $x$ that is given by 
\begin{itemize}
    \item $x_1 = 1$,
    \item $x_{i} = 2n-j+i-2$ for $2 \leq i \leq j+1$,
    \item $x_{i} = i$ for $j+2 \leq i \leq 2n-j-1$,
    \item $x_{i} = j+i-2n+2$ for $2n-j \leq i \leq 2n-1$, and
    \item $x_{2n} = 2n$
\end{itemize} belongs to $\mathsf{Pop}_{\mathrm{Weak}(B_n)}(\mathrm{Weak}(B_n))$. 
\end{lemma}
\begin{proof}
Construct $y$ by reversing each ascending run of $x$. Then $y$ is given by
\begin{itemize}
    \item $y_{i} = 2n-i$ for $1 \leq i \leq j$,
    \item $y_{j+1} = 1$,
    \item $y_{i} = 2n+1-i$ for $j+2 \leq i \leq 2n-j-1$,
    \item $y_{2n-j} = 2n$, and
    \item $y_{i} = 2n+2-i$ for $2n+1-j \leq i \leq 2n$. 
\end{itemize}
Then $\mathrm{asc}_1(y)$ starts at $y_1$ and end at $y_{j+1}$, while $\mathrm{asc}_2(y)$ starts at $y_{j+2}$ and ends at $y_{2n-j-1}$, and $\mathrm{asc}_3(y)$ starts at $y_{2n-j}$ and ends at $y_{2n}$. By reversing the descending runs of $y$, we verify that $\mathsf{Pop}_{\mathrm{Weak}(B_n)}(y) = x$ as sought. 
\end{proof}

Next, we show that any permutation $x \in \mathsf{Pop}_{\mathrm{Weak}(B_n)}(\mathrm{Weak}(B_n))$ that satisfies $x_1 = 1$ and $|\mathscr{U}_{\mathrm{Weak(B_n)}}(x)| = n-1$ either belongs to the class of permutations described in the statement of \Cref{one_one} or \Cref{one_two}.

\begin{lemma} \label{one_three}
Suppose that $x =1 \cdot (y+1) \cdot 2n \in \mathsf{Pop}_{\mathrm{Weak}(B_n)}(\mathrm{Weak(B_n)})$ and that $|\mathscr{U}_{\mathrm{Weak}(B_n)}(x)| = n-1$. If $y \not\in \mathsf{Pop}_{\mathrm{Weak}(B_{n-1})}(\mathrm{Weak}(B_{n-1}))$, then $x$ must be given by  
\begin{itemize}
    \item $x_1 = 1$,
    \item $x_{i} = 2n-j+i-2$ for $2 \leq i \leq j+1$,
    \item $x_{i} = i$ for $j+2 \leq i \leq 2n-j-1$,
    \item $x_{i} = j+i-2n+2$ for $2n-j \leq i \leq 2n-1$, and
    \item $x_{2n} = 2n$
\end{itemize} for some $1 \leq j \leq n-1$. 
\end{lemma}
\begin{proof}
    Suppose that $z \in \mathrm{Weak}(B_{n-1})$ is given by reversing each ascending run of $y$. Because $y \not \in \mathsf{Pop}_{\mathrm{Weak}(B_{n-1})}(\mathrm{Weak}(B_{n-1}))$, it must be that $\mathsf{Pop}_{\mathrm{Weak}(B_{n-1})}(z) \ne y$. In addition, because we assumed that $|\mathscr{U}_{\mathrm{Weak}(B_{n})}(x)| = n-1$, we have that $|\mathscr{U}_{\mathrm{Weak}(B_{n-1})}(y)| = n-2$. Thus, for $\mathsf{Pop}_{\mathrm{Weak}(B_{n-1})}(z) \ne y$, any number that appears in $\mathrm{asc}_1(y)$ must be greater than any number that appears in $\mathrm{asc}_{2}(y)$. Now because $|\mathscr{U}_{\mathrm{Weak}(B_{n-1})}(y)| = n-2$, any number that does not appear in $\mathrm{asc}_1(y)$ and $\mathrm{asc}_{-1}(y)$ must appear in $\mathrm{mid}(y)$. So for any number in $\mathrm{asc}_1(y)$ to be greater than any number in $\mathrm{asc}_2(y)$, the numbers that appear in $\mathrm{asc}_{1}(y)$ must be in the set $\{2n-j-1, 2n-j, \ldots, 2n-2\}$ for some $1 \leq j \leq n-1$. Because the numbers that appear in $\mathrm{asc}_{1}(y)$ appear in ascending order, we must have $y_{i} = 2n-3-j+i$ for all $1 \leq i \leq j$. Accordingly, the numbers that appear in $\mathrm{asc}_{-1}(y)$ are fixed to be $y_{i} = j+i-2n+3$ for all $2n-j-1 \leq i \leq 2n-2$. Now because $|\mathscr{U}_{\mathrm{Weak}(B_{n-1})}(y)| = n-2$, the numbers that do not appear in $\mathrm{asc}_{1}(y)$ and $\mathrm{asc}_{-1}(y)$ must appear in $\mathrm{mid}(y)$ in ascending order. Therefore, we have that $y_{i} = i$ for all $j+1 \leq i \leq 2n-j-2$. The statement of the lemma now follows. 
\end{proof}

We now proceed to the count the number of permutations $x \in \mathsf{Pop}_{\mathrm{Weak}(B_n)}(\mathrm{Weak}(B_n))$ that satisfy $|\mathscr{U}_{\mathrm{Weak}(B_{n})}(x)| = n-1$ and $2 \leq x_{1} \leq n$. We begin by proving a property that holds for all $x \in \mathrm{Weak}(B_n)$ such that $|\mathscr{U}_{\mathrm{Weak}(B_{n})}(x)| = n-1$. 

\begin{lemma} \label{two_one}
    Suppose that $x \in \mathrm{Weak}(B_n)$ satisfies $|\mathscr{U}_{\mathrm{Weak}(B_{n})}(x)|=n-1$. If some $j \geq n+1$ satisfies $\mathrm{ind}_{x}(j) \leq n$, then $j$ must appear in $\mathrm{asc}_{1}(x)$. 
\end{lemma} 
\begin{proof}
    Suppose otherwise. Then the index of the last number that appears in the ascending run that contains $j$ must be at most $n$, because $x_{2n+1-\mathrm{ind}_x(j)} = 2n+1-j$. Therefore, $|\mathscr{U}_{\mathrm{Weak}(B_{n})}(x)| \leq n-2$. We reach a contradiction.  
\end{proof}

We give an explicit formula for the number of permutations $x \in \mathsf{Pop}_{\mathrm{Weak}(B_n)}(\mathrm{Weak}(B_n))$ that satisfy $|\mathscr{U}_{\mathrm{Weak}(B_{n})}(x)| = n-1$ and $2 \leq x_{1} \leq n$. 

\begin{lemma} \label{two_count}
    The number of permutations $x \in \mathsf{Pop}_{\mathrm{Weak}(B_n)}(\mathrm{Weak}(B_n))$ that satisfy $x_{1} = i$ for some $2 \leq i \leq n$ and $|\mathscr{U}_{\mathrm{Weak}(B_{n})}(x)| = n-1$ is $2^{i-1} \cdot 3^{n-i}$.  
\end{lemma}
\begin{proof}
    Because $|\mathscr{U}_{\mathrm{Weak}(B_{n})}(x)| = n-1$, all entries of $x$ belong to one of $\mathrm{asc}_{1}(x)$, $\mathrm{asc}_{-1}(x)$, or $\mathrm{mid}(x)$. In addition, $\mathrm{asc}_{1}(x)$ must be distinct from $\mathrm{asc}_{-1}(x)$. Because $\mathrm{asc}_{1}(x)$, $\mathrm{mid}(x)$, and $\mathrm{asc}_{-1}(x)$ are ascending runs, after we fix the set of numbers that appear in each of $\mathrm{asc}_{1}(x)$, $\mathrm{mid}(x)$, and $\mathrm{asc}_{-1}(x)$, the permutation $x$ is fixed. Now because $x \in B_{n}$, the entries of $\mathrm{asc}_{-1}(x)$ must be fixed after the entries of $\mathrm{asc}_{1}(x)$ are fixed. In addition, every other number that does not appear in $\mathrm{asc}_{1}(x)$ or $\mathrm{asc}_{-1}(x)$ must appear in $\mathrm{mid}(x)$. Therefore, fixing the numbers that appear in $\mathrm{asc}_{1}(x)$ fixes $x$. 

    Since $x \in \mathrm{Weak}(B_n)$, for any $j \geq n+1$, either $\mathrm{ind}_{x}(j) \leq n$ or $ \mathrm{ind}_x(2n+1-j) \leq n$. First, suppose that $2n+1-j > x_{1}$. If $\mathrm{ind}_{x}(j) \leq n$, then $j$ must belong to $\mathrm{asc}_1(x)$ by \Cref{two_one}. On the other hand, if $\mathrm{ind}_x(2n+1-j) \leq n$, then either $2n+1-j$ belongs to $\mathrm{asc}_{1}(x)$ or $\mathrm{mid}(x)$.  
    
    Now suppose that $x_1 \geq 2n+1-j$. If $\mathrm{ind}_{x}(j) \leq n$, then $j$ must belong to $\mathrm{asc}_1(x)$ by \Cref{two_one}. But if $\mathrm{ind}_{x}(2n+1-j) \leq n$, then $2n+1-j$ must belong to $\mathrm{mid}(x)$, because $x_{1} > 2n+1-j$. 

     Collating, there are at most $2^{i-1} \cdot 3^{n-i}$ ways of fixing the numbers that belong to $\mathrm{asc}_{1}(x)$. Then $x$ is given by ordering the numbers that appear in each of $\mathrm{asc}_1(x)$, $\mathrm{mid}(x)$, and $\mathrm{asc}_{-1}(x)$ in ascending order and concatenating them. Now for each way of fixing the numbers that belong to $\mathrm{asc}_{1}(x)$, we can check that the numbers that we fixed to appear in $\mathrm{asc}_1(x)$, $\mathrm{asc}_{-1}(x)$, and $\mathrm{mid}(x)$ do appear in the appropriate ascending runs, because the largest number that we fixed to be in $\mathrm{asc}_1(x)$ is larger than the smallest number that we fixed to be in $\mathrm{asc}_2(x)$.
     
    Now for each $x$, we have that $x_{1}<(\mathrm{asc}_2(x))_{\ell(2)}$. Therefore, the preimage of $x$ with respect to $\mathsf{Pop}_{\mathrm{Weak}(B_n)}$ is given by reversing each ascending run of $x$. The statement of the lemma now follows. 
\end{proof}

Similarly, we provide a formula that counts the number of $x \in \mathsf{Pop}_{\mathrm{Weak}(B_n)}(\mathrm{Weak}(B_n))$ that satisfies $n+1 \leq x_{1} \leq 2n$ and $|\mathscr{U}_{\mathrm{Weak}(B_n)}(x)| = n-1$. 

\begin{lemma} \label{three_count}
     The number of permutations $x \in \mathsf{Pop}_{\mathrm{Weak}(B_n)}(\mathrm{Weak}(B_n))$ that satisfy $x_{1} = i$ for some $n+1 \leq i \leq 2n$ and $|\mathscr{U}_{\mathrm{Weak}(B_{n})}(x)| = n-1$ is $2^{2n-i}-1$.  
\end{lemma}
\begin{proof}
    As in the proof of \Cref{two_count}, because $|\mathscr{U}_{\mathrm{Weak}(B_{n})}(x)| = n-1$, fixing the entries that appear in $\mathrm{asc}_{1}(x)$ fixes the entries of $x$. Now every number greater than $x_{1}$ can either be in $\mathrm{asc}_{1}(x)$ or not. However, by \Cref{first}, $\ell(1) \ne 2n+1-x_1$. Therefore, there are at most $2^{2n-i}-1$ ways of fixing the entries of $\mathrm{asc}_{1}(x)$.

    For each way of fixing the numbers that appear in $\mathrm{asc}_1(x)$, the entries of $x$ are given by ordering the numbers that appear in each of $\mathrm{asc}_1(x)$, $\mathrm{mid}(x)$, and $\mathrm{asc}_{-1}(x)$ in ascending order and concatenating them. Now for each way of fixing the numbers that belong to $\mathrm{asc}_{1}(x)$, we can check that the numbers that we fixed to appear in $\mathrm{asc}_1(x)$, $\mathrm{asc}_{-1}(x)$, and $\mathrm{mid}(x)$ do appear in the appropriate ascending runs, because the largest number that we fixed to be in $\mathrm{asc}_1(x)$ is larger than the smallest number that we fixed to be in $\mathrm{asc}_2(x)$.
    
    Now for each $x$, we have that $x_{1}<(\mathrm{asc}_2(x))_{\ell(2)}$. Therefore, the preimage of $x$ with respect to $\mathsf{Pop}_{\mathrm{Weak}(B_n)}$ is given by reversing each ascending run of $x$. The statement of the lemma now follows. 
\end{proof}

We next proceed to the proof of \Cref{weak}. 

\begin{proof}[Proof of \Cref{weak}]
We induct on $n$. The base cases are clear. 

So we assume that the statement holds for $n-1$ and show that the statement holds for $n$. By \Cref{one_three} and the induction hypothesis, the number of $x \in \mathsf{Pop}_{\mathrm{Weak}(B_n)}(\mathrm{Weak}(B_n))$ that satisfies $x_1 = 1$ and $|\mathscr{U}_{\mathrm{Weak}(B_{n})}(x)|= n-1$ is given by $3^{n-1}-2(n-1)-1 + (n-1) = 3^{n-1}-n$. Next, by \Cref{two_count}, the number of $x \in \mathsf{Pop}_{\mathrm{Weak}(B_n)}(\mathrm{Weak}(B_n))$ that satisfies $2 \leq x_1 \leq n$ is given by $\sum_{i=2}^{n} 2^{i-1} \cdot 3^{n-i} = 2\cdot3^{n-1}-2^{n}$. Lastly, by \Cref{three_count}, the number of $x \in \mathsf{Pop}_{\mathrm{Weak}(B_n)}(\mathrm{Weak}(B_n))$ that satisfies $n+1 \leq x_1 \leq 2n$ is given by $\sum_{i= n+1}^{2n} (2^{2n-i}-1) = 2^{n}-n-1$. Collating, the number of $x \in \mathsf{Pop}_{\mathrm{Weak}(B_n)}(\mathrm{Weak}(B_n))$ that satisfies $|\mathscr{U}_{\mathrm{Weak}(B_{n})}(x)| = n-1$ is given by $(3^{n-1}-n)+(2\cdot3^{n-1}-2^{n})+(2^{n}-n-1) = 3^{n}-2n-1$ as sought. 
\end{proof}

\subsection{\texorpdfstring{$\mathrm{Tam}(B_n)$}{Tam(Bn)}} 

In this section, we settle a conjecture by Defant and Williams \cite{conjecture} on the coefficients of $\mathsf{Pop}(\mathrm{Tam}(B_n);q)$. We first introduce notations.

Suppose that $x = x_{1} x_{2} \cdots x_{2n} \in \mathrm{Tam}(B_{n})$. Let $\mathrm{half}(x)$ be the subsequence of $x$ that contains all the numbers that are at least $n+1$. For example, if $x = 65718243$, then $\mathrm{half}(x) = 6578$. Furthermore, let $\mathrm{half}_{k}(x) = x_{i_{k}(x)+1} \cdots x_{i_{k}(x) + m_{k}(x)}$ be the $k^{\text{th}}$ maximal substring from the left that contains numbers that are at least $n+1$. In our example $x = 65718243$, we have that $\mathrm{half}_{1}(x) = 657$ and $\mathrm{half}_{2}(x) = 8$. Let the \textit{complement $\mathrm{half}^{c}_{k}(x)$} of $\mathrm{half}_{k}(x)$ be the sequence $(2n+1- x_{i_{k}(x) + m_{k}(x)}) \cdot (2n+1-x_{i_{k}(x) +m_{k}(x)-1}) \cdots (2n+1-x_{i_{k}(x)+1})$. Continuing with our example of $x = 65718243$, we have that $\mathrm{half}^{c}_{1}(x) = 243$ and $\mathrm{half}^{c}_{2}(x) = 1$. Furthermore, let $\mathrm{len}_{k}(x) = \mathrm{len}(\mathrm{half}_{k}(x))$. For $x = 65718243$, we have that $\mathrm{len}_{1}(x) = 3$ and $\mathrm{len}_{2}(x) = 1$. 

 
We first establish a few necessary conditions on the permutations that lie in the image of $\mathsf{Pop}_{\mathrm{Tam}(B_n)}$. 

\begin{lemma} \label{nec1}
    Let $x \in \mathrm{Tam}(B_n)$. Then $\mathrm{ind}_{\mathsf{Pop}_{\mathrm{Tam}(B_{n})}(x)}(2n) \geq n+1$. 
\end{lemma}
\begin{proof}
    Suppose otherwise. Then $\mathrm{ind}_{x}(2n) \leq n$ as neither reversing descending runs nor swapping the positions of the numbers that correspond to $3$ and $1$ in a $\overline{31}2^{*}$ or a $2\overline{31}^{*}$ pattern can decrease the index of $2n$.

    Then we claim that $x_{i} \geq n+1$ for all $i$ such that $\mathrm{ind}_{x}(2n) < i \leq n$. Suppose otherwise. Then $x_{2n+1-i} \geq n+1$, because $x_{i}+x_{2n+1-i} = 2n+1$. Therefore, $2n$, $x_{i}$, and $x_{2n+1-i}$ form a $312^{*}$ pattern. But $x$ must avoid such a pattern as $x \in \mathrm{Tam}(B_n)$. Our claim thus follows. 

    Next we claim that $x_{i-1} > x_{i}$ for all $i$ such that $\mathrm{ind}_{x}(2n) < i \leq n$. Suppose otherwise. Then because $x_{i} \geq n+1$ by our previous claim, $2n,$ $x_{i-1},$ and $x_{i}$ form a $312^{*}$ pattern. But $x$ must avoid such a pattern as $x \in \mathrm{Tam}(B_n)$. The claim thus follows. 

    Now collating our two claims with the condition that $x_{i} + x_{2n+1-i} = 2n+1$, we have that $x_{\mathrm{ind}_{x}(2n)} > x_{\mathrm{ind}_{x}(2n)+1} > \cdots > x_{2n+1-\mathrm{ind}_{x}(2n)}$. Therefore, $x_{\mathrm{ind}_{x}(2n)}$ and $x_{2n+1-\mathrm{ind}_{x}(2n)}$ must belong to the same descending run of $x$, resulting in
    \begin{equation} \label{eq1b}
        \mathrm{ind}_{\mathrm{rev}(x)}(2n) \geq n+1.
    \end{equation}
     Now, because the sequence of swaps that remove a $\overline{31}2^{*}$ or a $2\overline{31}^{*}$ pattern can only increase the index of $2n$, we have 
    \begin{equation} \label{eq2b}
        \mathrm{ind}_{\mathsf{Pop}_{\mathrm{Tam}(B_{n})}(x)}(2n) \geq \mathrm{ind}_{\mathrm{rev}(x)}(2n).
    \end{equation}
    Combining \Cref{eq1b} and \Cref{eq2b}, we reach a contradiction as sought. 
\end{proof}

The following is another necessary condition on the permutations that lie in the image of $\mathsf{Pop}_{\mathrm{Tam}(B_n)}$. 

\begin{lemma} \label{nec2}
Let $x \in \mathrm{Tam}(B_n)$. Then for any $k$, we have that
\begin{equation}
    \mathrm{red}(\mathrm{half}_{k}(\mathsf{Pop}_{\mathrm{Tam}(B_n)}(x))) \in \mathsf{Pop}_{\mathrm{Tam}(A_{\mathrm{len}_{k}(x)-1})}(\mathrm{Tam}(A_{\mathrm{len}_{k}(x)-1}).
\end{equation} \end{lemma}

For the proof of \Cref{nec2}, we make use of the following theorem by Hong \cite{carina}. A \emph{double descent} of a permutation is a $\overline{321}$ pattern. 

\begin{theorem}[Hong, \cite{carina}]\label{hong1}
Let $x \in \mathrm{Tam}(A_{n-1})$. Then $\mathsf{Pop}_{\mathrm{Tam}(A_{n-1})}(x)$ ends with $n$ and does not contain a double descent. 
\end{theorem}

From \Cref{hong1}, we see that to prove \Cref{nec2}, we must show that the permutation $\mathrm{red}(\mathrm{half}_{k}(\mathsf{Pop}_{\mathrm{Tam}(B_n)}(x)))$ ends with $\mathrm{len}_{k}(x)$ and has no double descents. To do so, we prove the following series of auxiliary lemmas. We first show that for $x \in \mathrm{Tam}(B_n)$, the operations $\mathrm{half}$ and $\mathrm{rev}$ are commutative. 

\begin{lemma} \label{rev_half}
For $x \in \mathrm{Tam}(B_n)$, we have that $\mathrm{half}(\mathrm{rev}(x)) = \mathrm{rev}(\mathrm{half}(x))$. 
\end{lemma}
\begin{proof}
It suffices to show that $\mathrm{ind}_{x}(y_{i-1}) + 1 = \mathrm{ind}_{x}(y_{i})$ for every descent $y_{i-1} > y_{i}$ in $y = \mathrm{half}(x)$ as if so, then $y_{i}$ and $y_{i-1}$ will belong to the same descending run of $x$, from which the statement of the lemma follows. Suppose otherwise. Then because $y_{i-1}$ and $y_{i}$ are adjacent in $y$, we have that $x_{\mathrm{ind}_{x}(y_{i-1}) + 1} \leq n$. Therefore, $y_{i-1}$, $x_{\mathrm{ind}_{x}(y_{i-1}) + 1},$ and $y_{i}$ form a $312^{*}$ pattern. But $x$ must avoid such a pattern as $x \in \mathrm{Tam}(B_n)$. We reach a contradiction as sought. 
\end{proof}

We next prove the following constructive statement. 

\begin{lemma} \label{construction}
Suppose that $x \rhd_{\mathrm{Tam}(A_{n-1})} y$. Then for $z \in B_{n}$ such that $\mathrm{red}(\mathrm{half}(z)) = x$, there exists some $v^{0},v^1, v^2, \ldots, v^{\ell} \in B_{n}$ such that $z = v^{0} \rhd_{\mathrm{Tam}(B_n)} v^1 \rhd_{\mathrm{Tam}(B_n)}  \cdots \rhd_{\mathrm{Tam}(B_n)} v^{\ell}$ and $\mathrm{red}(\mathrm{half}(v^{\ell})) = y$. 
\end{lemma}
\begin{proof}
Because $x \rhd_{\mathrm{Tam}(A_{n-1})} y$, we begin by writing $x$ as $ x = X ca YbZ$ and $y$ as $y = X ac Y b Z$ for some sequence of numbers $X$, $Y,$ and $Z$ and numbers $a < b < c$. Now in $z$, the entries $c+n$, $z_{i}$, and $a+n$ form a $312^{*}$ pattern for all $\mathrm{ind}_{z}(c+n) + 1 \leq i \leq \mathrm{ind}_{z}(a+n) - 1$, because $c$ and $a$ are adjacent in $x$, which means that $z_{i} \leq n$. Successively define $v^{i}$ from $v^{i-1}$ by swapping the order of $c+n$ with the number that succeeds it until $c+n$ and $a+n$ are adjacent. Once $a+n$ and $c+n$ are adjacent, $a+n$, $c+n$, and $b+n$ form a $\overline{31}2^{*}$ pattern. Thus, we finish by swapping the order of $c+n$ and $a+n$. Our claim now follows. 
\end{proof}

\begin{example}
We demonstrate the statement of \Cref{construction} through the following example. Suppose that $x = 3142 \rhd_{\mathrm{Tam}(A_3)} 1342=y$. Take $z = 37145826$ such that $\mathrm{red}(\mathrm{half}(z)) = x$. Then by taking $z=37145826 \rhd_{\mathrm{Tam}(B_4)} 31745286  \rhd_{\mathrm{Tam}(B_4)} 31472586 \rhd_{\mathrm{Tam}(B_4)} 31427586 \rhd_{\mathrm{Tam}(B_4)} 31245786$, we arrive at $31245786$, which satisfies $\mathrm{red}(\mathrm{half} (31245786)) = 1342$ as sought.
\end{example}

We now repeatedly apply \Cref{construction} to arrive at the following corollary. 

\begin{corollary} \label{repeat}
 For $x \in \mathrm{Tam}(B_n)$, we have that 
 \begin{equation*}
     \mathrm{red}(\mathrm{half}(\pi_{\mathrm{Tam}(B_n)_{\downarrow}}(x))) = \pi_{\mathrm{Tam}(A_{n-1})_{\downarrow}}(\mathrm{red}(\mathrm{half}(x))).
 \end{equation*}
\end{corollary}
\begin{proof}
    From iteratively applying \Cref{construction}, we can find some $z \equiv_{\mathrm{Tam}(B_n)} x$ such that 
    \begin{equation} \label{eq3b}
        \mathrm{red}(\mathrm{half}(z)) = \pi_{\mathrm{Tam}(A_{n-1})_{\downarrow}}(\mathrm{red}(\mathrm{half}(x))). 
    \end{equation}
    Then because $\pi_{\mathrm{Tam}(A_{n-1})_{\downarrow}}(\mathrm{half}(x)) \in \mathrm{Tam}(A_{n-1})$, the permutation $z \in B_{n}$ avoids big $312^{*}$ patterns. Therefore, the relative ordering of any two numbers that are both at least $n+1$ remains invariant from $z$ to $\pi_{\mathrm{Tam}(B_n)_{\downarrow}}(z)$. Thus,  
    \begin{equation} \label{eq4b}
        \mathrm{half}(z) = \mathrm{half}(\pi_{\mathrm{Tam}(B_n)_{\downarrow}}(z)).
    \end{equation}
    
    Now combining \Cref{eq3b} and \Cref{eq4b}, we have that
    \begin{equation} \label{eq5b}
        \pi_{\mathrm{Tam}(A_{n-1})_{\downarrow}}(\mathrm{red}(\mathrm{half}(x)))= \mathrm{red}(\mathrm{half}(\pi_{\mathrm{Tam}(B_n)_{\downarrow}}(z))). 
    \end{equation}
    
    In addition, because $z \equiv_{\mathrm{Tam}(B_n)} x$, we have that
    \begin{equation} \label{eq6b}
        \mathrm{half}(\pi_{\mathrm{Tam}(B_n)_{\downarrow}}(z)) = \mathrm{half}(\pi_{\mathrm{Tam}(B_n)_{\downarrow}}(x)).
    \end{equation}
    
    Now combining \Cref{eq5b} and \Cref{eq6b}, we arrive at the statement of the corollary as desired. 
\end{proof}

\begin{example}
We demonstrate the statement of \Cref{repeat} through the following example. Suppose that $x= 37154826$. We have $\pi_{\mathrm{Tam}(A_3)}(\mathrm{half}(x)) = \pi_{\mathrm{Tam}(A_3)}(3142) = 1342$. Following the steps of \Cref{construction}, we can find $x= 37154826 \rhd_{\mathrm{Tam}(B_4)} 31754286 \rhd_{\mathrm{Tam}(B_4)} 31572486 =z$. Now, because $\mathrm{red}(\mathrm{half}(z)) = 1342 \in \mathrm{Tam}(A_3)$, we have that $z$ avoids big $312^{*}$ patterns. We must now perform a sequence of swaps to $z$ to remove all small $312^{*}$ patterns: $z = 31572486  \rhd_{\mathrm{Tam}(B_4)} 31527486 \rhd_{\mathrm{Tam}(B_4)} 31527486 \rhd_{\mathrm{Tam}(B_4)} 31254786 \rhd_{\mathrm{Tam}(B_4)} 31245786 = \pi_{\mathrm{Tam}(B_n)_{\downarrow}}(x)$. We check that $\mathrm{red}(\mathrm{half}(\pi_{\mathrm{Tam}(B_4)_{\downarrow}}(x))) = \mathrm{red}(\mathrm{half}(31245786)) = 1342 = \pi_{\mathrm{Tam}(A_3)_{\downarrow}}(3142) = \pi_{\mathrm{Tam}(A_3)_{\downarrow}}(\mathrm{red}(\mathrm{half}(x)))$ as sought. 
\end{example}

We next establish a lemma on the relative order of the numbers that appear in $\mathrm{half}_{k}(x)$ across $k$ for $x \in \mathrm{Tam}(B_n)$.

\begin{lemma} \label{prop}
    For $x \in \mathrm{Tam}(B_n)$, the smallest number in $\mathrm{half}_{k}(x)$ is greater than the largest number in $\mathrm{half}_{k-1}(x)$.
\end{lemma}
\begin{proof}
Suppose otherwise. Let $c$ be the largest number in $\mathrm{half}_{k-1}(x)$. Next, let $a$ be a number appearing strictly between $\mathrm{half}_{k-1}(x)$ and $\mathrm{half}_{k}(x)$. By definition, $a \leq n$. Now, let $b$ be the smallest number in $\mathrm{half}_{k}(x)$. Then $c,$ $a,$ and $b$ form a $312^*$ pattern in $x$. Therefore, we reach a contradiction as sought.  
\end{proof}

In a related direction, we establish a lemma on the relative order of the numbers in $y$ for $y \in \mathrm{rev}(\mathrm{Tam}(A_{n-1}))$. 

\begin{lemma} \label{order}
For $y \in \mathrm{Tam}(A_{n-1})$ and $x= \mathrm{rev}(y)$, we have that $ y_{n}-1 < x_{j}$ for all $\mathrm{ind}_{x}(y_{n}-1) < j$.  
\end{lemma}
\begin{proof}
Suppose otherwise. Because $x_{j} \leq y_{n}-1$ and $\mathrm{ind}_x(y_n-1) < j$, the numbers $x_{j}$ and $y_{n}-1$ must have belonged to different descending runs of $y$. Thus, the relative order of $x_{j}$ and $y_{n}-1$ must have been preserved by $\mathrm{rev}$. Therefore, we have that $\mathrm{ind}_{y}(y_{n}-1) < \mathrm{ind}_{y}(x_{j})$. 

Now for all $k < \mathrm{ind}_{y}(x_{j})$, we have that $y_{k} < y_{n}$, as otherwise, $y_{k},$ $x_{j}$, and $y_{n}$ form a $312$ pattern in $y$. Therefore, $y_{n}-1$ is the largest number that appears in $y_{1}y_{2} \cdots y_{\mathrm{ind}_{y}(x_{j})}$. Now suppose that $y_{n}-1$ does not belong to the last descending run of $y_{1}y_{2} \cdots y_{\mathrm{ind}_{y}(x_{j})}$. Let $y_{\ell}$ be the end of the descending run in $y_{1}y_{2} \cdots y_{\mathrm{ind}_{y}(x_{j})}$ that contains $y_{n}-1$. Because $y_{n}-1$ is the maximal number among $y_{1}, y_{2}, \ldots, y_{\mathrm{ind}_{y}(x_{j})},$ we have that $y_{\ell} \ne y_{n}-1$. Then $y_{n}-1, y_{\ell},$ and $y_{\ell+1}$ form a $312$ pattern in $y$. Therefore, $y_{n}-1$ must belong to the last descending run of $y_{1}y_{2} \cdots y_{\mathrm{ind}_{y}(x_{j})}$. But if so, then $y_{n}-1$ and $x_{j}$ must belong to the same descending run of $y$. Thus, it cannot be that $\mathrm{ind}_{x}(y_{n}-1) < j$. As a result, we reach a contradiction. 
\end{proof}

We next establish a lemma on the relative ordering of the numbers in $\mathrm{rev}(y)$ for $y \in \mathrm{Tam}(B_n)$. 

\begin{lemma} \label{rev_order}
    For $x= \mathrm{rev}(y)$ such that $y \in \mathrm{Tam}(B_n)$, if $x_{i} \geq n+1$ and $x_{i+1} \leq n$, then $x_{i} > x_{j}$ for all $j < i$.  
\end{lemma}
\begin{proof}
    Suppose otherwise. We first show that $x_{i+1} \ne 0$. Let $x_{\ell}$ be the smallest number in the same descending run as $2n$. Then $2n$ must belong to the rightmost descending run of $\mathrm{half}(y)$ as otherwise, $2n$, $x_{\ell}$, and $x_{\ell+1}$ form a $312^{*}$ pattern. Thus, $\mathrm{half}(x)_{n} = 2n$ as desired. 
    
    Now, because $x_{i} \geq n+1$ and $x_{i+1} \leq n$, it must be that $x_{i}$ is the largest number in its descending run in $y$. Therefore, $x_{i}$ and $x_{j}$ must belong to different descending runs in $y$. Because reversing descending runs does not change the relative ordering of the numbers that belong to different descending runs, we have that $\mathrm{ind}_{y}(x_{j}) < \mathrm{ind}_{y}(x_{i})$ from $\mathrm{ind}_{x}(x_j) = j < i = \mathrm{ind}_{x}(x_i)$.
    
     Lastly, because $y$ must avoid $312^{*}$ patterns, there cannot be any number less than $x_{i}$ between $x_{j}$ and $x_{i}$ in $y$. Therefore, $x_{i}$ cannot be the largest element in its descending run of $y$. We thus reach a contradiction. 
\end{proof}

We next use \Cref{rev_order} to establish the analog of \Cref{rev_order} on $x \in \mathsf{Pop}_{\mathrm{Tam}(B_n)}(\mathrm{Tam}(B_n))$. 

\begin{lemma} \label{analog}
    For $x= \mathsf{Pop}_{\mathrm{Tam}(B_n)}(y)$ such that $y \in \mathrm{Tam}(B_n)$, if $x_{i} \geq n+1$ and $x_{i+1} \leq n$, then $x_{i} > x_{j}$ for all $j < i$. 
\end{lemma}
\begin{proof}
    Consider the sequence of swaps from $\mathrm{rev}(y)$ to $x =\mathsf{Pop}_{\mathrm{Tam}(B_n)}(y)$. Because the order of the swaps does not matter, we first remove all small $312^{*}$ patterns. We begin by observing that for all $k$ such that $\mathrm{half}_{k}(y)$ does not contain $y_{2n}$, the numbers that belong to the last descending run of $\mathrm{half}_{k}(y)$ appear in increasing order in $\mathrm{half}_{h(k)}(\mathrm{rev}(y))$, while any other number in $\mathrm{half}_{k}(y)$ appears in $\mathrm{half}_{h(k)-1}(\mathrm{rev}(y))$ for some $h(k)$. Combining this observation with \Cref{prop}, the small $312^{*}$ patterns in $\mathrm{rev}(y)$ have the element that corresponds to $3$ belong to $\mathrm{half}_{k}(y)$ but not in the last descending run of $\mathrm{half}_{k}(y)$ and the element that corresponds to $2$ belong to the last descending run of $\mathrm{half}_{k}(y)$ for some $k$. Furthermore, from \Cref{prop}, we see that the smallest number in $\mathrm{half}_{h(k)}(\mathrm{rev}(y))$ is $\mathrm{half}_{h(k)}(\mathrm{rev}(y))_1$. 

    Now, in $\mathrm{half}_{\ell(k)-1}(\mathrm{rev}(y))$, if a number greater than $\mathrm{half}_{\ell(k)}(\mathrm{rev}(y))_1$ appears to the left of a number less than $\mathrm{half}_{\ell(k)}(\mathrm{rev}(y))_1$, then the two numbers with $\mathrm{half}_{\ell(k)}(\mathrm{rev}(y))_1$ form a $312^{*}$ pattern. Therefore, we can swap the order between the numbers in $\mathrm{half}_{\ell(k)-1}(\mathrm{rev}(y))$ such that any number smaller than $\mathrm{half}_{\ell(k)}(\mathrm{rev}(y))_1$ appears to the left of any number greater than $\mathrm{half}_{\ell(k)}(\mathrm{rev}(y))_1$.

    Next, any number greater than $\mathrm{half}_{h(k)}(\mathrm{rev}(y))_1$ in $\mathrm{half}_{h(k)-1}(\mathrm{rev}(y))$, any number in between $\mathrm{half}_{h(k)-1}(\mathrm{rev}(y))$ and $\mathrm{half}_{h(k)}(\mathrm{rev}(y))$, and $\mathrm{half}_{h(k)}(\mathrm{rev}(y))_1$ form a $312^{*}$ pattern. Therefore, we can perform a sequence of swaps such that the numbers that are greater than $\mathrm{half}_{h(k)}(\mathrm{rev}(y))_1$ that used to belong to $\mathrm{half}_{h(k)-1}(\mathrm{rev}(y))$ are now just in front of $\mathrm{half}_{h(k)}(\mathrm{rev}(y))_1$. 

    Let the permutation that results after performing the swaps described above for all $k$ be $\mathrm{interim}(y)$. Because any number in $\mathrm{half}_{k-1}(\mathrm{interim}(y))$ is smaller than any number in $\mathrm{half}_{k}(\mathrm{interim}(y))$, there are no small $312^{*}$ patterns in $\mathrm{interim}(y)$. Moreover, because of \Cref{order} and \Cref{rev_order}, we see that in $\mathrm{interim}(y)$, the last element in $\mathrm{half}_{k}(\mathrm{interim}(y))$ for any $k$ is greater than any number that appears to the left of it. 
    
    From $\mathrm{interim}(y)$ to $\mathsf{Pop}_{\mathrm{Tam}(B_n)}(y)$, we must rid of all big $312^{*}$ patterns. To do so, all we have to do is switch the order of elements that belong to the same $\mathrm{half}_{k}(\mathrm{interim}(y))$. Now, because the last element in $\mathrm{half}_{k}(\mathrm{interim}(y))$ is the greatest element in $\mathrm{half}_{k}(\mathrm{interim}(y))$, its position cannot change. The statement of the lemma now follows. 
\end{proof}

\begin{example}

We demonstrate the argument presented in the proof of \Cref{analog} in action through the following example. Let $$y = 7 \cdot 1 \cdot 10 \cdot 11 \cdot 9 \cdot 8 \cdot 5\cdot 4 \cdot 2 \cdot 3 \cdot 12 \cdot 6.$$ Then $$\mathrm{rev}(y) = 1 \cdot 7 \cdot 10 \cdot 2 \cdot 4 \cdot 5 \cdot 8 \cdot 9 \cdot 11 \cdot 3 \cdot 6 \cdot 12.$$
We now perform a sequence of swaps that remove all the small $312^{*}$ patterns. Because $10$ in $\mathrm{half}_1(\mathrm{rev}(y))$ is greater than $8$ and $9$ in $\mathrm{half}_2(\mathrm{rev}(y))$, we swap the position of $10$ with its subsequent number until we have $$ \mathrm{interim}(y) = 1 \cdot 7 \cdot 2 \cdot 4 \cdot 5 \cdot 3 \cdot 10 \cdot 8 \cdot 9 \cdot 11 \cdot 6 \cdot 12.$$
Finally, we remove big $312^{*}$ patterns by switching the order of the numbers in $\mathrm{half}_2(\mathrm{interim}(y))$. We arrive at 
$$\mathsf{Pop}_{\mathrm{Tam}(B_6)}(y) =1 \cdot 7 \cdot 2 \cdot 4 \cdot 3 \cdot 5 \cdot 8 \cdot 10 \cdot 9 \cdot 11 \cdot 6 \cdot 12 .$$
As in the statement of \Cref{analog}, we see that the rightmost numbers in $\mathrm{half}_{1}(\mathsf{Pop}_{\mathrm{Tam}(B_6)}(y))$, $\mathrm{half}_{2}(\mathsf{Pop}_{\mathrm{Tam}(B_6)}(y))$, and $\mathrm{half}_{3}(\mathsf{Pop}_{\mathrm{Tam}(B_6)}(y))$ are larger than any number to its left. 
\end{example}

We now proceed to the proof of \Cref{nec2}. 

\begin{proof}[Proof of \Cref{nec2}]

Because of \Cref{hong1} and \Cref{repeat}, we see that the permutation $\mathrm{red}(\mathrm{half}_{k}(\mathsf{Pop}_{\mathrm{Tam}(B_n)}(x)))$ cannot contain double descents. Moreover, from \Cref{analog}, we see that $\mathrm{red}(\mathrm{half}_{k}(\mathsf{Pop}_{\mathrm{Tam}(B_n)}(x)))$ must end with $\mathrm{len}_{k}(x)$. Thus, the statement of the lemma follows from \Cref{hong1}. 
\end{proof}

\Cref{nec1} and \Cref{nec2} give two necessary conditions on $y \in \mathsf{Pop}_{\mathrm{Tam}(B_n)}(\mathrm{Tam}(B_n))$. We next show that these two conditions are also sufficient. 

\begin{lemma} \label{suff}
      Suppose that $x \in \mathrm{Tam}(B_n)$. Then $x \in \mathsf{Pop}_{\mathrm{Tam}(B_n)}(\mathrm{Tam}(B_n))$ if $\mathrm{ind}_{y}(2n) \geq n+1$ and $\mathrm{red}(\mathrm{half}_{k}(\mathsf{Pop}_{\mathrm{Tam}(B_n)}(x))) \in \mathsf{Pop}_{\mathrm{Tamari(A_{\mathrm{len}_{k}(x)-1})}}(\mathrm{Tamari}(A_{\mathrm{len}_{k}(x)-1})$ for all $k$.   
\end{lemma}

Before proceeding to the proof of \Cref{suff}, we establish the following auxiliary lemma. 

\begin{lemma}\label{end1}
Suppose that $x \in \mathsf{Pop}_{\mathrm{Tam}(A_{n-1})}(\mathrm{Tam}(A_{n-1}))$. Then there exists $y \in \mathrm{Tam}(A_{n-1})$ such that $y_{n} = 1$ and $\mathsf{Pop}_{\mathrm{Tam}(A_{n-1})}(y) = x$. 
\end{lemma} 
\begin{proof}
    We induct on $n$. The base case is clear. Let $\mathrm{ind}_{x}(1) = k$. Then because $x$ cannot contain a $312$ pattern, we have that $\{x_{1}, x_{2}, \ldots, x_{k-1}\} = \{2, 3, \ldots, k\}$ and $\{x_{k+1}, x_{k+2}, \ldots, x_{n}\} = \{k+1, k+2, \ldots, n\}$. Now, from \Cref{hong1}, $x$ cannot contain a double descent. Moreover, $x$ must avoid $312$ patterns. Thus, $x_{k-1} = k$ as otherwise, we either introduce a $312$ or a $\overline{321}$ pattern. In addition, by \Cref{hong1}, we know that $x_{n} = n$. 
    
    Therefore, it follows from \Cref{hong1} that $\mathrm{red}(x_{1}x_{2} \cdots x_{k-1}) \in \mathsf{Pop}_{\mathrm{Tam}(A_{k-2})}(\mathrm{Tam}(A_{k-2}))$. Let $z \in \mathrm{Tam}(A_{k-2})$ such that $\mathsf{Pop}_{\mathrm{Tam}(A_{k-2})}(z) = \mathrm{red}(x_{1}x_{2} \cdots x_{k-1})$.  Similarly, we have that $\mathrm{red}(x_{k+1}x_{k+2} \cdots x_{n}) \in \mathsf{Pop}_{\mathrm{Tam}(A_{n-k-1})}(\mathrm{Tam}(A_{n-k-1}))$. From the induction hypothesis, let $w \in \mathrm{Tam}(A_{n-k-1})$ such that $\mathsf{Pop}_{\mathrm{Tam}(A_{n-k-1})}(w) = \mathrm{red}(x_{k+1}x_{k+2} \cdots x_{n})$ and $w_{n-k} = 1$. 

    Now, consider $(z+1) \cdot (w+k) \cdot 1$. We first observe that $(z+1) \cdot (w+k) \cdot 1 \in \mathrm{Tam}(A_{n-1})$, because $z \in \mathrm{Tam}(A_{k-2})$, $w \in \mathrm{Tam}(A_{n-k-1})$, and any number in $z+1$ is less than any number in $w+k$. Thus, we can apply $\mathsf{Pop}_{\mathrm{Tam}(A_{n-1})}$ to $(z+1) \cdot (w+k) \cdot 1$. 

    In $\mathrm{rev}((z+1) \cdot (w+k) \cdot 1)$, we see that $1$ appears to the left of $k+1$. Therefore, we start by switching the order of $1$ and any number at least $k+2$ that appears to the left of it, because the two numbers and $k+1$ form a $312$ pattern. Now, after the sequence of swaps, all numbers in a $312$ pattern either belong to $z+1$ or $w+k$. Therefore, we see that 
    \begin{align}
        \mathsf{Pop}_{\mathrm{Tam}(A_{n-1})}((z+1) \cdot (w+k) \cdot 1) & = (\mathsf{Pop}_{\mathrm{Tam}(A_{n-1})}(z)+1) \cdot 1 \cdot (\mathsf{Pop}_{\mathrm{Tam}(A_{n-1})}(w)+k)   \nonumber  \\ & = x_{1}x_{2} \cdots x_{k-1} \cdot 1 \cdot x_{k+1} x_{k+2} \cdots x_{n} \\ & = x. \nonumber
    \end{align}
The statement of the lemma now follows. 
\end{proof}

\begin{example}
    We demonstrate the statement of \Cref{end1} through the following example. Let $x =24351768$. Now, because $\mathsf{Pop}_{\mathrm{Tam}(A_3)}(\mathrm{red}(4532)) = 1324 = \mathrm{red}(2435)$, $\mathsf{Pop}_{\mathrm{Tam}(A_2)}(\mathrm{red}(786)) = 213 = \mathrm{red}(768)$, and $\mathrm{ind}_{\mathrm{red}(786)}(1) = \mathrm{len}(786)$, following the proof of \Cref{end1}, we set $y = 45327861$.  Then, 
    \begin{align*}
        \mathsf{Pop}_{\mathrm{Tam}_{A_7}}(45327861) &= \pi_{\mathrm{Tam}(A_7)}(42357168)  \\ &= \pi_{\mathrm{Tam}(A_7)}(42351768) \\ & = \pi_{\mathrm{Tam}(A_3)}(4235) \cdot 1 \cdot \pi_{\mathrm{Tam}(A_2)}(768)  \\ &= (\mathsf{Pop}_{\mathrm{Tam}_{A_3}}(3421) + 1) \cdot 1 \cdot (\mathsf{Pop}_{\mathrm{Tam}_{A_2}}(213)+5)  \\ &= 24351768,
    \end{align*} as desired. 
\end{example}

We now proceed to the proof of \Cref{suff}. 

\begin{proof}[Proof of \Cref{suff}]
We begin by writing $x$ as: 
\begin{equation*}
    x = \mathrm{half}^{c}_{\ell}(x) \cdot \mathrm{half}_{1}(x) \cdot \mathrm{half}^{c}_{\ell-1}(x) \cdot \mathrm{half}_{2}(x) \cdots \mathrm{half}_{\ell}(x) \cdot \mathrm{half}^{c}_{1}(x),
\end{equation*}
where $\mathrm{half}_{\ell}(x)$ and $\mathrm{half}^{c}_{\ell}(x)$ are possibly empty. 

Now, by \Cref{end1}, there exists
\begin{equation*}
    y =  \mathrm{half}_{1}(y) \cdot \mathrm{half}^{c}_{\ell}(y) \cdot \mathrm{half}_{2}(y) \cdot \mathrm{half}^{c}_{\ell-1}(y) \cdots \mathrm{half}_{\ell}(y) \cdot \mathrm{half}^{c}_{1}(y),
\end{equation*}
such that for each $k$,
\begin{itemize}
    \item the set of numbers that appear in $\mathrm{half}_{k}(y)$ coincide with the set of numbers that appear in $\mathrm{half}_{k}(y)$,
    \item $\mathsf{Pop}_{\mathrm{Tam}(A_{\mathrm{len}_{y}(k) - 1})}(\mathrm{red}(\mathrm{half}_{k}(y))) = \mathrm{red}(\mathrm{half}_{k}(x))$, and
    \item $\mathrm{ind}_{\mathrm{red}(\mathrm{half}_{k}(y))}(1)= \mathrm{len}_{y}(k)$. 
\end{itemize}
It follows from \Cref{prop} that $y \in \mathrm{Tam}(B_n)$. Thus, we can apply $\mathsf{Pop}_{\mathrm{Tam}(B_n)}$ to $y$. 

We first reverse the descending runs of $y$. From \Cref{prop}, we know that in the small $312^{*}$ patterns that appear in $\mathrm{rev}(y)$, the element that corresponds to $3$ appears in $\mathrm{half}_{k}(y)$ but not in its last descending run, the element that corresponds to $1$ appears in $\mathrm{half}^{c}_{\ell+1-k}(y)$, and the element that corresponds to $2$ appears in the last descending run of $\mathrm{half}_{k}(y)$ for some $k$. 

Recall that we let the subsequence $\mathrm{half}_{k}(y)$ of $y$ correspond to $y_{i(k)+1}$ $\cdot$ $y_{i(k)+2}$ $\cdots$ $y_{i_{k}(y) + \mathrm{len}_k(y)}$. Then $\mathrm{half}^{c}_{\ell+1-k}(y)$ is given by $y_{i_{k}(y)+\mathrm{len}_k(y)+1}$ $\cdot$ $y_{i_{k}(y)+\mathrm{len}_k(y) + 2}$ $\cdots$ $y_{i_{k}(y) + \mathrm{len}_k(y) + \mathrm{len}_{\ell+1-k}(y)}$. We see from \Cref{prop} that even if $\mathrm{half}^{c}_{\ell}(x)$ is empty, the descending runs that begin in $\mathrm{half}_{\ell+1-k}(x)$ end in $\mathrm{half}^{c}_{k}(x)$ without running over to $\mathrm{half}_{\ell+2-k}(x)$. Thus, we have that $\{y_{i_{k}(y)+1}, \ldots, y_{i_{k}(y) + \mathrm{len}_k(y)+\mathrm{len}_{\ell+1-k}(y)}\} = \{\mathrm{rev}(y)_{i_{k}(y)+1}, \ldots, \mathrm{rev}(y)_{i_{k}(y) + \mathrm{len}_{k}(y)+\mathrm{len}_{\ell+1-k}(y)}\}$. 

Now, write
\begin{equation*}
    \mathrm{rev}(y)_{i_{k+1}(y)} \cdot \mathrm{rev}(y)_{{i_{k}(y)}+2} \cdots \mathrm{rev}(y)_{i_{k}(y) + \mathrm{len}_k(y)+\mathrm{len}_{\ell+1-k}(y)} = s \cdot t \cdot u \cdot v
\end{equation*}
such that 
\begin{itemize}
    \item the numbers that appear in $s$ are the numbers that appear in $\mathrm{half}_{k}(y)$ but not in its last descending run,
    \item the numbers that appear in $t$ are the numbers that appear in the first descending run of $\mathrm{half}^{c}_{\ell+1-k}(x)$,  
    \item the numbers that appear in $u$ are the numbers that appear in the last descending run of $\mathrm{half}_{k}(y)$, and
    \item the numbers that appear in $v$ are the numbers that appear in $\mathrm{half}^{c}_{\ell+1-k}(x)$ but not in its first descending run.   
\end{itemize}

Then it follows that 
\begin{equation*}
    \mathrm{rev}(y)_{i_{\ell+1-k}+1} \cdot \mathrm{rev}(y)_{i_{\ell+1-})+2} \cdots \mathrm{rev}(y)_{i_{\ell+1-k}(y) + \mathrm{len}_{k}(y)+\mathrm{len}_{\ell+1-k}} = v^{c} \cdot u^{c} \cdot t^{c} \cdot s^{c}. 
\end{equation*}

Now, because $\mathrm{ind}_{\mathrm{red}(\mathrm{half}_{k}(y))}(1)= \mathrm{len}_{k}(y)$, we have that $u_{1}$ is the smallest number that appears in $\mathrm{half}_{k}(y)$. Thus, any number in $s$, any number in $t$, and $u_{1}$ form a $312^{*}$ pattern. Therefore, we can switch the order of elements in $t$ and $s$ such that after the sequence of swaps, we have $s$, $t$, $u$, and $v$ appear in the order of $t \cdot s \cdot u \cdot v$. Consequently, $s^{c}$, $t^{c}$, $u^{c}$, and $v^{c}$ appear in the order of $v^{c} \cdot u^{c} \cdot s^{c} \cdot t^{c}$.

Also, because $\mathrm{ind}_{\mathrm{red}(\mathrm{half}_{\ell+1-k}(x))}(1) = \mathrm{len}_{\ell+1-k}(x)$, we see that $t_{\mathrm{len}(t)}$ is the smallest number contained in $\mathrm{half}^{c}_{\ell+1-k}(x)$. So, $t_{\mathrm{len}(t)}$, any number in $\mathrm{half}_{k}(x)$, and any number in $v$ form a $231^{*}$ pattern. Thus, we can switch the order of the numbers in $\mathrm{half}_{k}(x)$ and any number in $v$ such that after the swaps, we have $s$, $t$, $u$, $v$ appear in the order of $t \cdot v \cdot s \cdot u$. Equivalently, $t \cdot v \cdot s \cdot u =\mathrm{rev}(\mathrm{half}^{c}_{\ell+1-k}(y)) \cdot  \mathrm{rev}(\mathrm{half}_{k}(y))$. Similarly, we have that $s^{c}$, $t^{c}$, $u^{c}$, and $v^{c}$ appear in the order of $u^{c} \cdot s^{c} \cdot v^{c} \cdot t^{c} = \mathrm{rev}(\mathrm{half}^{c}_{k}(y)) \cdot \mathrm{rev}(\mathrm{half}_{\ell+1-k}(y))$. 

Let the permutation at the end of these series of swaps be $\mathrm{aux}(y)$. Then we see that 
\begin{equation*}
    \mathrm{aux}(y) = \mathrm{rev}(\mathrm{half}^{c}_{\ell}(y)) \cdot \mathrm{rev}(\mathrm{half}_{1}(y)) \cdots \mathrm{rev}(\mathrm{half}^{c}_{1}(y)) \cdot \mathrm{rev}(\mathrm{half}_{\ell}(y)). 
\end{equation*}

Then because $\mathsf{Pop}_{\mathrm{Tam}(A_{\mathrm{len}_{k}(x) - 1})}(\mathrm{red}(\mathrm{half}_{k}(y))) = \mathrm{red}(\mathrm{half}_{k}(x))$, 
we have that 
\begin{align*}
    \pi_{\mathrm{Tam}(B_n)\downarrow}(\mathrm{aux}(y)) &=  \mathrm{half}^{c}_{\ell}(x) \cdot \mathrm{half}_{1}(x) \cdot \mathrm{half}^{c}_{\ell-1}(x) \cdot \mathrm{half}_{2}(x) \cdots \mathrm{half}_{\ell}(x) \cdot \mathrm{half}^{c}_{1}(x) \\ & = x.
\end{align*}
The statement of the lemma now follows. 
\end{proof}

\begin{example}

We demonstrate the statement of \Cref{suff} in action through the following example. Let $$x = 1 \cdot 7 \cdot 2 \cdot 4 \cdot 3 \cdot 5 \cdot 8 \cdot 10 \cdot 9 \cdot 11 \cdot 6 \cdot 12.$$ Then we see that
$$y = 7 \cdot 1 \cdot 10 \cdot 11 \cdot 9 \cdot 8 \cdot 5\cdot 4 \cdot 2 \cdot 3 \cdot 12 \cdot 6$$ satisfies the conditions that $y$ must satisfy in the proof of \Cref{suff}. For example, we show that the conditions hold for $k =2$: 
\begin{itemize}
    \item $\mathrm{half}_2(x)$ and $\mathrm{half}_2(y)$ both consist the numbers $8$, $9$, $10$, and $11$,
    \item $\mathsf{Pop}_{\mathrm{Tam}(A_3)}(\mathrm{red}(\mathrm{half}_{2}(y))) = \mathsf{Pop}_{\mathrm{Tam}(A_3)}(\mathrm{red}(10 \cdot 11 \cdot 9 \cdot 8)) = \mathsf{Pop}_{\mathrm{Tam}(A_3)}(3 4  2 1)= \pi_{\mathrm{Tam}(A_3) \downarrow}(3124) = 1324 =  \mathrm{red}(8 \cdot 10 \cdot 9 \cdot 11) = \mathrm{red}(\mathrm{half}_{2}(x))$, and
    \item $\mathrm{ind}_{\mathrm{red}(\mathrm{half}_{2}(y))}(1)= \mathrm{ind}_{\mathrm{red}(10 \cdot 11 \cdot 9 \cdot 8)}(1)= \mathrm{ind}_{(3421)}(1)= 4 = \mathrm{len}_{2}(y)$.
\end{itemize}

By following the steps outlined in the proof, we see that $$ \mathrm{aux}(y) = 1 \cdot 7 \cdot 2 \cdot 4 \cdot 5 \cdot 3 \cdot 10 \cdot 8 \cdot 9 \cdot 11 \cdot 6 \cdot 12.$$ 

Finally, we see that $$\mathsf{Pop}_{\mathrm{Tam}(B_6)}(y) = 1 \cdot 7 \cdot 2 \cdot 4 \cdot 3 \cdot 5 \cdot 8 \cdot 10 \cdot 9 \cdot 11 \cdot 6 \cdot 12 = x$$ as sought. 

\end{example}

Lastly, we combine \Cref{nec1}, \Cref{nec2}, and \Cref{suff} to arrive at the necessary and sufficient conditions for an element $x \in \mathrm{Tam}(B_n)$ to be in $\mathsf{Pop}_{\mathrm{Tam}(B_n)}(\mathrm{Tam}(B_n))$. 

\begin{theorem} \label{character}
For $x \in \mathrm{Tam}(B_n)$, we know that $x \in \mathsf{Pop}_{\mathrm{Tam}(B_n)}(\mathrm{Tam}(B_n))$ if and only if $\mathrm{ind}_{x}(2n) \geq n+1$ and $\mathrm{red}(\mathrm{half}_{k}(\mathsf{Pop}_{\mathrm{Tam}(B_n)}(x))) \in \mathsf{Pop}_{\mathrm{Tam}(A_{\mathrm{len}_{k}(x)-1})}(\mathrm{Tam}(A_{\mathrm{len}_{k}(x)-1}))$ for all $k$. 
\end{theorem}

We now proceed to the proof of \Cref{tamaribn}, in which we make use of the following result by Hong \cite{carina}. 

\begin{theorem}[Hong, \cite{carina}] \label{ancount}
For all $n \geq 1$, 
\begin{equation*}
    \mathsf{Pop}(\mathrm{Tam}(A_n);q) = \sum_{k=0}^{n} \frac{1}{k+1} \binom{2k}{k} \binom{n}{2k} q^{n-k}. 
\end{equation*}
\end{theorem}

\begin{proof}[Proof of \Cref{tamaribn}]
We begin by setting $\mathcal{M}$ to be the set of permutations $z$ that belong in $\mathsf{Pop}_{\mathrm{Tam}(A_{\mathrm{len}(z)-1})}(\mathrm{Tam}(A_{\mathrm{len}(z)-1}))$. In addition, let
\begin{equation}
    M(x, y) := \sum_{z \in \mathcal{M}} x^{\mathrm{len}(z)}y^{\lvert \mathscr{U}_{\mathrm{Tam}(A_{\mathrm{len}(z)-1})(z)} \rvert}.
\end{equation}

From \Cref{ancount} and the OEIS Sequence A055151 \cite{oeis}, we have that 
\begin{equation} \label{mexp}
    M(x,y) = x\left( - \frac{-1+x + \sqrt{1-2x+x^2-4x^2y^2}}{2x^2y^2} \right). 
\end{equation}

Next, we let $\mathcal{N}$ be the set of $z$ that belong to $\mathsf{Pop}_{\mathrm{Tam}(B_{\mathrm{len}(z)/2})}(\mathrm{Tam}(B_{\mathrm{len}(z)/2}))$. As an analog of $M(x,y)$, we let  
\begin{equation}
    N(x, y) := \sum_{z \in \mathcal{N}} x^{\mathrm{len}(z)/2}y^{|\mathrm{des}(z)|}.
\end{equation}

Furthermore, let $\mathcal{P}$ be the set of $z$ that belong to $\mathsf{Pop}_{\mathrm{Tam}(B_{\mathrm{len}(z)})}(\mathrm{Tam}(B_{\mathrm{len}(z)}))$ and satisfy $z_1 \geq \mathrm{len}(z)/2 + 1$. We define
\begin{equation}
    P(x, y) := \sum_{z \in \mathcal{P}} x^{\mathrm{len}(z)/2}y^{|\mathrm{des}(z)|}.
\end{equation}

It then follows from \Cref{prop} and \Cref{character} that 
\begin{equation} \label{peq}
        P(x, y) = y^2M(x,y^2)^2 + y^3M(x,y^2) + \dots = \frac{y^2M(x,y^2)^2}{1-yM(x,y^2)}.
\end{equation}

In the same vein, we let $\mathcal{Q}$ be the set of $z$ that belong to $\mathsf{Pop}_{\mathrm{Tam}(B_{\mathrm{len}(z)})}(\mathrm{Tam}(B_{\mathrm{len}(z)}))$ and satisfy $z_1 \leq \mathrm{len}(z)/2$. Let 
\begin{equation}
    Q(x, y) := \sum_{z \in \mathcal{Q}} x^{\mathrm{len}(z)/2}y^{|\mathrm{des}(z)|}.
\end{equation}

It then follows from \Cref{prop} and \Cref{character} that
\begin{equation} \label{qeq}
    Q(x, y) = M(x,y^2) + yM(x,y^2)^2 + \dots = \frac{M(x,y^2)}{1-yM(x,y^2)}. 
\end{equation}

Now from \Cref{mexp}, \Cref{peq} and \Cref{qeq}, we have that 

\begin{align}
    N(x,y) &= P(x, y) + Q(x,y) \nonumber \\ 
    &= \frac{y^2M(x,y^2)^2}{1-yM(x,y^2)} + \frac{M(x,y^2)}{1-yM(x,y^2)}\\
    &= - \frac{-1+x+2x^2y^2 + \sqrt{1-2x + x^2 - 4x^2y^2}}{xy(-1+x+2xy + \sqrt{1-2x+x^2-4x^2y^2})}. \nonumber
\end{align}

Now if we let $L(x,y) = \left( 1- \left( \frac{2xy}{1-x}\right)^2 \right)^{\frac{1}{2}}$, then $N(x,y) = xN_0(x,y) + xy N_1(x,y)$, where 
\begin{equation}
    N_0(x,y) = \frac{1}{1-x} + \frac{1+x}{2x(1-x)} \left( \frac{1}{L(x,y)} -1 \right)
\end{equation}

and 

\begin{equation}
    N_1(x,y) = (1-x)^{-2} \left( \frac{x+1}{L(x,y)} - \frac{1}{2} \left( \frac{1-x}{xy} \right)^2 \left(1 - L(x,y) \right) \right). 
\end{equation}

We now simplify the expression of $N_0(x,y)$ by expanding $L(x,y)$ out of the denominator:  \begin{align} \label{simplify}
    N_0(x,y) &= \frac{1}{1-x} + \frac{1+x}{2x(1-x)} \left( \frac{1}{L(x,y)} -1 \right) \nonumber \\ &= \frac{1}{1-x} + \frac{1+x}{2x(1-x)} \sum_{k=1}^{\infty} \binom{2k}{k} \left( \frac{xy}{1-x} \right)^{2k} \\
    &= \frac{1}{1-x} + \frac{1}{2} \sum_{k= 1}^{\infty} \binom{2k}{k} (1-x)^{-2k-1} (1+x) x^{2k-1} y^{2k} \nonumber. 
\end{align}
Then because for $k \geq 1$, we have that
\begin{equation}
    \sum_{n =0}^{\infty} (n+1-k) \frac{n!}{(n+1-2k)!} x^n = \frac{1}{2}(2k)! (1-x)^{-2k-1} (1+x) x^{2k-1},
\end{equation} we can further simplify the expression of $N_0(x,y)$ from \Cref{simplify} as follows. 
\begin{align}
    N_0(x,y) &= \frac{1}{1-x} + \frac{1}{2} \sum_{k= 1}^{\infty} \binom{2k}{k} (1-x)^{-2k-1} (1+x) x^{2k-1} y^{2k} \nonumber \\ &= \sum_{n = 0}^{\infty} x^n + \sum_{k=1}^{\infty} \frac{y^{2k}}{(k!)^2} \sum_{n =0}^\infty (n+1-k) \frac{n!}{(n+1-2k)!} x^n  \\
    &= \sum_{n =0}^{\infty} \sum_{k = 0}^{\infty} (n+1-k) \frac{n!}{(k!)^2 (n+1-2k)!} x^n y^{2k} \nonumber \\
    &= \sum_{n =0}^{\infty} \sum_{k = 0}^{\infty} \binom{n}{k} \binom{n+1-k}{k} x^n y^{2k} \nonumber.
\end{align}

Similarly, by expanding out $L(x,y)$, we have that
\begin{align} \label{neq}
    N_1(x,y) & = (1-x)^{-2} \left( \frac{x+1}{L(x,y)} - \frac{1}{2} \left( \frac{1-x}{xy} \right)^2 \left(1 - L(x,y) \right) \right)  \\
    &= (1-x)^{-2} \left( (x+1) \sum_{k = 0}^\infty \binom{2k}{k} \left( \frac{xy}{1-x} \right)^{2k} - \sum_{k=0}^\infty \frac{1}{k+1} \binom{2k}{k} \left(\frac{xy}{1-x} \right)^{2k}\right) \nonumber \\
    &= (1-x)^{-2} \sum_{k = 0}^{\infty} \binom{2k}{k} \frac{1}{k+1} \left( (k+1)(x+1) -1 \right) \left( \frac{xy}{1-x} \right)^{2k} \nonumber.
\end{align}

Furthermore, because 
\begin{equation}
    \sum_{n =0}^{\infty} \binom{n}{2k} (n-k) x^n = (1-x)^{-2} \left( \frac{x}{1-x} \right)^{2k} \left( (k+1)x +k \right),
\end{equation} we can simplify the expression of $N_1(x,y)$ from \Cref{neq} further as follows. 
\begin{align}
    N_1(x,y) &= (1-x)^{-2} \sum_{k = 0}^{\infty} \binom{2k}{k} \frac{1}{k+1} \left( (k+1)(x+1) -1 \right) \left( \frac{xy}{1-x} \right)^{2k} \nonumber \\ &= \sum_{n = 0}^{\infty} \sum_{k=0}^{\infty} \frac{n!}{(k+1)! k! (n-2k)!} (n-k) x^n y^{2k} \\
    &= \sum_{n = 0}^{\infty} \sum_{k= 0}^{\infty} \binom{n}{k+1} \binom{n-k}{k} x^n y^{2k} \nonumber. 
\end{align}

Then because $N(x,y)= x N_0(x,y) + xy N_1(x,y)$, we have that 
\begin{equation} \label{final}
    N(x,y) = \sum_{n =0}^{\infty} \sum_{k = 0}^{\infty} \binom{n-1}{k} \binom{n-k}{k} x^n y^{2k} + \sum_{n = 0}^{\infty} \sum_{k= 0}^{\infty} \binom{n-1}{k} \binom{n-k}{k-1} x^n y^{2k-1}. 
\end{equation}

Now define
\begin{equation}
    K(x,y) := \sum_{z \in \mathcal{N}} x^{\mathrm{len}(z)/2}y^{|\mathscr{U}_{\mathrm{Tam}(B_{\mathrm{len}(z)/2)}}(z)|}. 
\end{equation}

Because for $z \in B_n$, we have that $|\mathscr{U}_{\mathrm{Tam}(B_n)}(z)| = n-k$ if and only if $\mathrm{des}(z) = 2k$ or $\mathrm{des}(z) = 2k-1$. Therefore, we have that 
\begin{align} \label{finall}
    K(x,y) &= \sum_{n=0}^{\infty}  \sum_{k=0}^{\infty} \left( \binom{n-1}{k} \binom{n-k}{k} + \binom{n-1}{k} \binom{n-k}{k-1} \right) x^{n}y^{n-k}\\
    &=  \sum_{n=0}^{\infty}  \sum_{k=0}^{\infty} \binom{n-1}{k} \binom{n-k+1}{k} x^{n}y^{n-k} \nonumber.
\end{align}
The statement of the theorem now follows from \Cref{finall}. 
\end{proof}

\subsection{\texorpdfstring{$J(\Phi^{+}_{A_n})$}{J(Phi+An)}} 

Here, we settle a conjecture on the coefficients of $\mathsf{Pop}(J(\Phi^{+}_{A_n});q)$ presented by Defant and Williams \cite{conjecture}. The special case $q=1$ was previously shown to be true by Sapounakis, Tasoulas, and Tsikouras \cite{greek}. We begin by citing a result of Sapounakis, Tasoulas, and Tsikouras \cite{greek} that characterizes $\mu \in J(\Phi^{+}_{A_n})$ that belongs to $\mathsf{Pop}^{\uparrow}_{J(\Phi^{+}_{A_n})}(J(\Phi^{+}_{A_n}))$. Recall that we call $(1,1)$ a $\emph{rise}$ and denote it by $\mathit{r}$ and $(1, -1)$ a \emph{fall} and denote it by $\mathit{f}$.
 
\begin{lemma}[Sapounakis, Tasoulas, and Tsikouras \cite{greek}] \label{aaaa} 
The path $\mu \in J(\Phi^{+}_{A_{n}})$ is in $\mathsf{Pop}^{\uparrow}_{J(\Phi^{+}_{A_n})}(J(\Phi^{+}_{A_n}))$ if and only if there are no consecutive steps of $\mathit{ffrr}$, and the only lattice points on $\mu$ on the $x-$axis are $(0, 0)$ and $(2n, 0)$.
\end{lemma}

We now prove \Cref{jayan} by refining the generating function presented in Sapounakis, Tasoulas, and Tsikouras \cite{greek}. 

\begin{proof}[Proof of \Cref{jayan}]
We start by defining $\mathcal{F}$ to be the set of Dyck paths $\mu$ that satisfy $\mu \in \mathsf{Pop}^{\uparrow}_{J(\Phi^{+}_{A_{s(\mu)}})}(J(\Phi^{+}_{A_{s(\mu)}}))$. We define
\begin{equation}
    F(x, y) := \sum_{\mu \in \mathcal{F}} x^{s(\mu)} y^{|\mathscr{D}_{J(\Phi^{+}_{A_{s(\mu)}})}(\mu)|}.
\end{equation}

Next, let $\mathcal{G}$ be the set of Dyck paths of any length that avoid consecutive steps of $ffrr$. Let
\begin{equation}
    G(x,y) := \sum_{\mu \in \mathcal{G}} x^{s(\mu)} y^{|\mathscr{D}_{J(\Phi^{+}_{A_{s(\mu)}})}(\mu)|}. 
\end{equation}

Then by \Cref{aaaa}, every nonempty $\mu \in \mathcal{F}$ can be written uniquely as $\textit{r} w \textit{f}$, where $w \in \mathcal{G}$. In addition, because the number of peaks of $w$ and $\mu$ are equal except for when $w$ is empty, we have that 
\begin{equation} \label{eq1}
    F(x,y) = 1 + x(G(x,y)-1) + xy.
\end{equation}

Furthermore, because every path $\mu$ can be uniquely written as $\mu = \textit{rf}v$, $\mu=\textit{r} w \textit{f}$, or $\mu = \textit{r}w\textit{frf}v$ for $v, w(\ne \varepsilon) \in \mathcal{G}$, we have that
\begin{equation} \label{eq2}
    G(x,y) = 1 + xy G(x,y) + x (G(x,y)-1) + x^2 y G(x, y) (G(x, y)-1).
\end{equation}

Now let 
\begin{equation} \label{eq3}
    H(x, y) := G(x, y) - 1. 
\end{equation} By solving for $H(x,y)$ from \Cref{eq1}, \Cref{eq2}, and \Cref{eq3}, we have that 
\begin{equation}
    H(x, y) = y\left(\frac{x(H(x,y)+1) + x^2 H(x, y)(H(x,y) + 1)}{1-x}\right).
\end{equation}
So from the Lagrange inversion theorem as seen in \cite{lagrange}, it follows that 
\begin{equation}
    [y^{k}] H(x,y) = \frac{1}{k} [u^{k-1}] \phi(u)^{k}, 
\end{equation}where $\phi(u) = \frac{x}{1-x} (u+1)(1+xu)$. 

Since
\begin{equation}
    [u^{k-1}] \phi(u)^{k}  = \sum_{j = 0}^{k} \binom{k}{j+1} \binom{k}{j}  \frac{x^{k+j}}{(1-x)^{k}},
\end{equation}
we have that 
\begin{equation}
    [x^{n}y^{k}] H(x,y) = \sum_{j = 0}^{n-k} \frac{1}{k} \binom{k}{j+1} \binom{k}{j} \binom{n-j-1}{n-k-j}.
\end{equation}
 It follows that 
 \begin{equation}
     H(x, y) = \sum_{n=0}^{\infty} \sum_{k=0}^{n} \sum_{j = 0}^{n-k} \frac{1}{k} \binom{k}{j+1} \binom{k}{j} \binom{n-j-1}{n-k-j} x^{n} y^{k}.
 \end{equation}

Now the statement of the theorem that 
\begin{equation}
    F(x, y) = \sum_{n=0}^{\infty} \sum_{k=0}^{n} \sum_{j=0}^{n-k+1} \frac{1}{k+1} \binom{k+1}{j} \binom{k+1}{j-1} \binom{n-j+1}{n-k-j+1}  x^{n}y^{k+1}
\end{equation}
follows from \Cref{eq1} and \Cref{eq3}.
\end{proof}

\subsection{\texorpdfstring{$J(\Phi^{+}_{B_n})$}{J(Phi+Bn)}} 

Here, we settle a conjecture by Defant and Williams \cite{conjecture} on the coefficients of $\mathsf{Pop}(J(\Phi^{+}_{B_{n}}); q)$. We characterize the paths $\mu \in J(\Phi^{+}_{B_{n}})$ that belong to $\mathsf{Pop}^{\uparrow}_{J(\Phi^{+}_{B_{n}})}(J(\Phi^{+}_{B_{n}}))$. 

\begin{lemma} \label{bbbb}
    The path $\mu \in J(\Phi^{+}_{B_{n}})$ is in $\mathsf{Pop}^{\uparrow}_{J(\Phi^{+}_{B_n})}(J(\Phi^{+}_{B_{n}}))$ if and only if there are no consecutive steps of $ffrr$, and the only lattice points on $\mu$ on the $x-$axis are $(0,0)$ and $(4n, 0)$. 
\end{lemma}
\begin{proof}
    For any $\nu \in J(\Phi^{+}_{B_{n}})$, we have that 
    \begin{equation} \label{eq4}
        \mathsf{Pop}^{\uparrow}_{J(\Phi^{+}_{B_{n}})}(\nu) = \mathsf{Pop}^{\uparrow}_{J(\Phi^{+}_{A_{2n}})}(\nu). 
    \end{equation}
    Because $J(\Phi^{+}_{B_{n}})$ is a sublattice of $J(\Phi^{+}_{A_{2n}})$, that any $\mu \in \mathsf{Pop}^{\uparrow}_{J(\Phi^{+}_{B_{n}})}(J(\Phi^{+}_{B_{n}}))$ must avoid consecutive steps of $ffrr$ and not touch the $x-$axis besides at $(0, 0)$ and $(4n,0)$ follows from \Cref{aaaa}.

    Conversely, assume that $\mu$ avoids $ffrr$ and does not intersect the $x-$axis besides at the two endpoints of the path. Let $\nu$ be the path given by transforming each peak in $\mu$ into a valley. It is clear that $\nu \in J(\Phi^{+}_{B_n})$, because $\mu$ is symmetric about $x=2n$ and does not intersect the $x-$axis besides at the two endpoints of the path. Now it is shown in the proof of \Cref{aaaa} in \cite{greek} that 
    \begin{equation}
        \mathsf{Pop}^{\uparrow}_{J(\Phi^{+}_{A_{2n}})}(\nu) =  \mu,
    \end{equation} because $\mu$ avoids $ffrr$. 
    It now follows from \Cref{eq4} that $\mu \in \mathsf{Pop}^{\uparrow}_{J(\Phi^{+}_{B_{n}})}(J(\Phi^{+}_{B_{n}}))$ as sought. 
\end{proof}

We next proceed to the proof of \Cref{jaybn}. We take a similar approach to our approach in \Cref{jayan}. 

\begin{proof}[Proof of \Cref{jaybn}]

As in the proof of \Cref{jayan}, we let $\mathcal{G}$ be the set of Dyck paths that avoid consecutive steps of \textit{ffrr} and let 
\begin{equation} 
    G(x, y) := \sum_{\mu \in \mathcal{G}} x^{s(\mu)} y^{|\mathscr{D}_{J(\Phi^{+}_{A_n})}(\mu)|}. 
\end{equation}

By solving for $G(x,y)$ explicitly from \Cref{eq2}, we have that 
\begin{equation} \label{eq5}
    G(x,y) = -\frac{-1+x+xy-x^2y + \sqrt{(-1+x)(4x^2y + (-1+x)(-1+xy)^2)}}{2x^2y}.
\end{equation}

Now let $\mathcal{I}$ be the set of Dyck paths $\mu$ that avoid $ffrr$ and touch the $x-$axis only at its endpoints. Let
\begin{equation}
    I(x, y) = \sum_{\mu \in \mathcal{I}} x^{s(\mu)}y^{p(\mu)}. 
\end{equation}

Because every path $\mu$ in $\mathcal{I}$ can be written uniquely as either $\mu = rfvrf$, $\mu = rvf$, $\mu = rwfrfvrfrwr$, or $\mu = rwfrfrwr$, where $v \in \mathcal{I}$ and $w (\ne \varepsilon) \in \mathcal{G}$, we have that 
\begin{equation} \label{eq6}
    I(x,y) = 1 + x^2 y I(x,y) + x I(x,y) + x^4 y I(x,y)(G(x^2,y)-1) + x^3 y (G(x^2, y)-1) -x + xy.
\end{equation}

Thus, solving for $I(x,y)$ from \Cref{eq5} and \Cref{eq6}, we have that 
\begin{equation} \label{eq6c}
    I(x,y) = \frac{1+2x-3x^2 +x^2y-x^4y-\sqrt{(-1+x^2)(4x^4y+(-1+x^2)(-1+x^2y)^2)}}{x(1-2x+x^2-x^2y+x^4y+\sqrt{(-1+x^2)(4x^4y+(-1+x^2)(-1+x^2y)^2)}}.
\end{equation}

Next let $\mathcal{J}$ be the set of Dyck paths $\mu$ that belong to $\mathsf{Pop}^{\uparrow}_{J(\Phi^{+}_{B_{s(\mu)/2}})}(J(\Phi^{+}_{B_{s(\mu)/2}}))$. Let 
\begin{equation}
    J(x,y) := \sum_{\mu \in \mathcal{J}} x^{s(\mu)} y^{|\mathscr{D}_{J(\Phi^{+}_{B_{s(\mu)/2}})}(\mu)|}.
\end{equation}

Now by \Cref{bbbb}, $\mu \in \mathcal{J}$ can be written as $r w f$, where $w \in \mathcal{I}$. Because $s(\mu)$ is even, $s(w)$ must be odd. Therefore, we have that 
\begin{equation} \label{eq7}
    J(x,y) = 1+ \frac{1}{2} \sqrt{x}(I(\sqrt{x},y) + I(-\sqrt{x},y)).
\end{equation}

From \Cref{eq6c} and \Cref{eq7}, we have that 

\begin{align} \label{eq8}
    J(x,y) =1 -\frac{\splitfrac{(x(-2+2x-y-2xy+3x^2y+xy^2-2x^2y^2+x^3y^2}{+ (2-y+xy)\sqrt{(x-1)(4x^2y+(x-1)(xy-1)^2)})}}{\splitfrac{(1-2x+x^2-2xy+2x^3y+x^2y^2-2x^3y^2+x^4y^2}{+ (1+x-xy+x^2y)\sqrt{(x-1)(4x^2y+(x-1)(xy-1)^2)})}}.
\end{align}

Now we set $z := \frac{xy}{x-1}$ and $w := xy+1$. Then \Cref{eq8} reduces to 
\begin{equation}
    J(x,y) = (w^2 + 4z)^{-\frac{1}{2}}. 
\end{equation}
From the binomial theorem, we have that
\begin{align*}
    J(x, y) & = \sum_{j=0}^{\infty} \binom{-\frac{1}{2}}{j} w^{-1-2j}(4z)^{j} \\ & = \sum_{j=0}^{\infty} \binom{-\frac{1}{2}}{j} (xy+1)^{-1-2j} (\frac{4xy}{x-1})^{j} \\ & = \sum_{j=0}^{\infty} \sum_{l=0}^{\infty} \sum_{m=0}^{\infty} \binom{-\frac{1}{2}}{j} \binom{-1-2j}{l} \binom{-j}{m} (-4xy)^{j}(xy)^{l}(-x)^{m} \\ & = \sum_{j=0}^{\infty} \sum_{l=0}^{\infty} \sum_{m=0}^{\infty} \binom{-\frac{1}{2}}{j} \binom{-1-2j}{l} \binom{-j}{m} (-4)^{j}(-1)^{m}(x)^{j+l+m}(y)^{j+l} \\ &= \sum_{n=0}^{\infty} \sum_{k = 0}^{\infty} \sum_{j=1}^{m} \binom{2j}{j} \binom{k+j}{k-j} \binom{n-k+j-1}{n-k} (-1)^{k-j} x^{n} y^{k}. 
\end{align*}

The statement of the theorem now follows. 

\end{proof}


\section{Future Directions}
\label{futuredirections}
In this paper, we resolved the conjectures proposed by Defant and Williams \cite{conjecture} concerning $\mathsf{Pop}(M; q)$ when $M$ is $\mathrm{Weak}(B_n)$, $\mathrm{Tam}(B_n)$, $J(\Phi^{+}(A_n))$, or $J(\Phi^{+}(B_n))$. Previously, their conjecture concerning $\mathsf{Pop}(M; q)$ for $M = \mathrm{Tam}(A_n)$ was proven by Hong \cite{carina}. The last of their six conjectures remains open; it concerns a formula for $\mathsf{Pop}(M; q)$ when $M=\mathrm{Camb}_{\mathrm{bi}}(A_n)$ is the Cambrian lattice of type $A_n$ arising from a bipartite Coxeter element. One can define this lattice to be the sublattice of $\mathrm{Weak(A_n)}$ consisting of the permutations $x_1\cdots x_{n+1}\in S_{n+1}$ such that 
\begin{itemize}
\item there do not exist indices $i_1<i_2<i_3$ with $x_{i_3}$ odd and $x_{i_2}<x_{i_3}<x_{i_1}$;
\item there do not exist indices $i_1<i_2<i_3$ with $x_{i_1}$ even and $x_{i_3}<x_{i_1}<x_{i_2}$.
\end{itemize}

\begin{conjecture}[Defant and Williams \cite{conjecture}]\label{cambrian}
$$\mathsf{Pop}(\mathrm{Camb}_{\mathrm{bi}}(A_n); q) = \sum_{k =1}^{\lfloor \frac{n+3}{2} \rfloor} \frac{k (-1)^{k-1}}{n-k+3} \sum_{j=0}^{n-k+3}  \binom{j}{n-j+3} \binom{n-k+3}{j} q^{j-2} .$$
\end{conjecture}

We end with a conjecture that gives a characterization of the permutations in $\mathrm{Camb}_{\mathrm{bi}}(A_n)$ that belong to $\mathsf{Pop}_{\mathrm{Camb}_{\mathrm{bi}}(A_n)}(\mathrm{Camb}_{\mathrm{bi}}(A_n))$. 

\begin{conjecture} \label{char_conj}
    Suppose that $x \in \mathrm{Camb}_{\mathrm{bi}}(A_n)$. Then $x \in \mathrm{Pop}_{\mathrm{Camb}_{\mathrm{bi}(A_n)}}(\mathrm{Camb}_{\mathrm{bi}}(A_n))$ if and only if $x$ satisfies all of the following: 
    \begin{itemize}
        \item $x$ avoids the pattern $\overline{321}$.
        \item $\mathrm{ind}_{x}(2k) < \mathrm{ind}_{x}(2k+1)$ for any $k$ such that $3 \leq 2k+1 \leq n$. 
        \item $\mathrm{ind}_{x}(2k) < \mathrm{ind}_{x}(2k+3)$ for any $k$ such that $5 \leq 2k+3 \leq n$.
        \item $\mathrm{ind}_x(n-2) < \mathrm{ind}_x(n)$ if $n \equiv 0 \pmod{2}$. 
        \item $\mathrm{ind}_{x}(1) < \mathrm{ind}_{x}(3)$. 
    \end{itemize}
\end{conjecture}

\section*{Acknowledgements}
\noindent This research was conducted at the 2022 University of Minnesota Duluth REU, supported by Jane Street Capital, the NSA (grant number H98230-22-1-0015), and the NSF (grant number DMS-2052036). In addition, Choi was supported by the Harvard Herchel Smith Undergraduate Science Research Program, and Sun was supported by the Harvard College Research Program. The authors are indebted to Joe Gallian for his dedication and organizing the University of Minnesota Duluth REU. We thank Amanda Burcroff, Jonas Iskander, and Colin Defant for their invaluable feedback and advice on this paper.

\nocite{*}
\printbibliography

\end{document}